\documentclass{article}

\usepackage{fullpage}
\usepackage{authblk}
\usepackage{natbib}
\usepackage{amsthm}
\usepackage{amsmath}
\usepackage{amsfonts}
\usepackage{amssymb}
\usepackage{mathtools}
\usepackage{thmtools}

\declaretheorem[name=Lemma,numberwithin=section]{lemma}
\declaretheorem[name=Theorem,sibling=lemma]{theorem}
\declaretheorem[name=Proposition,sibling=lemma]{proposition}

\declaretheorem[style=definition,name=Definition,sibling=lemma]{definition}

\declaretheorem[style=remark,name=Remark,sibling=lemma]{remark}
\declaretheorem[style=remark,name=Example,sibling=lemma]{example}

\usepackage{xparse}
\makeatletter
\NewDocumentEnvironment{savetheorem}{m +b}{
  \begin{theorem}\label{#1}#2\end{theorem}%
  \expandafter\gdef\csname rest@body@#1\endcsname{#2}%
}{}%
\newcommand{\printtheorem}[1]{%
  \par\noindent\textbf{Theorem~\ref{#1}. }\itshape
  \csname rest@body@#1\endcsname\par\normalfont
}

\NewDocumentEnvironment{saveprop}{m +b}{%
  \begin{proposition}\label{#1}#2\end{proposition}%
  \expandafter\gdef\csname rest@body@#1\endcsname{#2}%
}{}%
\newcommand{\printprop}[1]{%
  \par\noindent\textbf{Proposition~\ref{#1}. }\itshape
  \csname rest@body@#1\endcsname\par\normalfont
}
\makeatother

\usepackage{fix-cm}
\usepackage{booktabs}
\usepackage{lmodern}
\usepackage[mathscr]{eucal}
\usepackage{dirtytalk}
\usepackage[shortlabels]{enumitem}
\usepackage[T1]{fontenc}
\usepackage{tikz}
\usetikzlibrary{shapes.geometric,arrows,positioning}
\usepackage{multirow}
\usepackage{footmisc}
\usepackage{wrapfig}
\allowdisplaybreaks
\usepackage[colorlinks=true,linkcolor=red,citecolor=blue]{hyperref}
\usepackage{cleveref}
\newcommand{\sE}[0]{\mathscr{E}}
\newcommand{\sF}[0]{\mathscr{F}}
\newcommand{\sH}[0]{\mathscr{H}}
\newcommand{\sG}[0]{\mathscr{G}}
\newcommand{\sP}[0]{\mathscr{P}}

\definecolor{fannyComment}{RGB}{150, 10, 250}

\newcommand{\EE}[0]{\mathbb{E}}
\newcommand{\KK}[0]{\mathbb{K}}

\newcommand{\PP}[0]{\mathbb{P}}
\newcommand{\QQ}[0]{\mathbb{Q}}

\newcommand{\cC}[0]{\mathcal{C}}
\newcommand{\cD}[0]{\mathcal{D}}
\newcommand{\cM}[0]{\mathcal{M}}
\newcommand{\cU}[0]{\mathcal{U}}

\newcommand{\bC}[0]{\mathbf{C}}
\newcommand{\bF}[0]{\mathbf{F}}
\newcommand{\bN}[0]{\mathbf{N}}
\newcommand{\bS}[0]{\mathbf{S}}
\newcommand{\bT}[0]{\mathbf{T}}
\newcommand{\bU}[0]{\mathbf{U}}
\newcommand{\bV}[0]{\mathbf{V}}
\newcommand{\bW}[0]{\mathbf{W}}
\newcommand{\bX}[0]{\mathbf{X}}
\newcommand{\bY}[0]{\mathbf{Y}}
\newcommand{\bZ}[0]{\mathbf{Z}}
\newcommand{\bPA}[0]{\mathbf{PA}}

\newcommand{\bc}[0]{\boldsymbol{c}}
\newcommand{\bn}[0]{\boldsymbol{n}}
\newcommand{\bs}[0]{\boldsymbol{s}}
\newcommand{\bt}[0]{\boldsymbol{t}}
\newcommand{\bu}[0]{\boldsymbol{u}}

\newcommand{\bw}[0]{\boldsymbol{w}}
\newcommand{\bx}[0]{\boldsymbol{x}}
\newcommand{\bz}[0]{\boldsymbol{z}}

\newcommand{\iPP}[2]{\PP^{\textnormal{do}(#1,#2)}}
\newcommand{\iEE}[2]{\EE^{\textnormal{do}(#1,#2)}}
\newcommand{\iKK}[2]{\KK^{\textnormal{do}(#1,#2)}}
\newcommand{\iK}[3]{K^{\textnormal{do}(#1,#2)}_{#3}}

\newcommand{\fact}[1]{#1^\textnormal{F}}
\newcommand{\cfact}[1]{#1^\textnormal{CF}}

\newcommand{\ind}[0]{\mathbf{1}}

\DeclareRobustCommand{\indep}[1]{\mathpalette\IndepAux{#1}}
\newcommand{\IndepAux}[2]{
  \ifx#1\displaystyle \def\k{-10mu}\else
  \ifx#1\textstyle    \def\k{-9mu}\else
  \ifx#1\scriptstyle  \def\k{-7mu}\else
                       \def\k{-6mu}\fi\fi\fi
  \nonscript\mskip\thickmuskip\nobreak
  \mathord{\perp\mkern\k\perp}_{\!#2}%
  \nobreak\nonscript\mskip\thickmuskip
}
\newcommand{\asequal}[1]{\stackrel{#1}{=}}

\newcommand{\exam}[0]{\text{Exam}}
\newcommand{\class}[0]{\text{Class}}
\newcommand{\sky}[0]{\text{Sky}}
\newcommand{\starr}[0]{\text{Star}}

\title{Counterfactual spaces}

\author[1]{Junhyung Park}
\author[1]{Fanny Yang}
\author[2]{Thomas Icard}
\affil[1]{ETH Z\"urich}
\affil[2]{Stanford University}
\date{}

\begin{document}
\maketitle
\begin{abstract}
We mathematically axiomatise the stochastics of \emph{counterfactuals}, by introducing two related frameworks, called \emph{counterfactual probability spaces} and \emph{counterfactual causal spaces}, which we collectively term \emph{counterfactual spaces}. They are, respectively, probability and causal spaces whose underlying measurable spaces are products of world-specific measurable spaces. In contrast to more familiar accounts of counterfactuals founded on causal models, we do not view interventions as a necessary component of a theory of counterfactuals. As an alternative to Pearl's celebrated \say{ladder of causation}, we view counterfactuals and interventions are orthogonal concepts, respectively mathematised in counterfactual probability spaces and causal spaces. The two concepts are then combined to form counterfactual causal spaces. At the heart of our theory is the notion of shared information between the worlds, encoded completely within the probability measure and causal kernels, and whose extremes are characterised by \emph{independence} and \emph{synchronisation} of worlds. Compared to existing frameworks, counterfactual spaces enable the mathematical treatment of a strictly broader spectrum of counterfactuals.
\end{abstract}

\section{Introduction}\label{sec:introduction}
Counterfactual thinking is central to human cognition, behaviour and actions. Accordingly, it has received much attention within various academic disciplines. The tradition of possible worlds semantics has long shaped the philosophical discussions of counterfactuals \citep{goodman1947problem,lewis1973counterfactuals,lewis1986plurality,stalnaker2003ways}, while psychologists have studied how imagining counterfactual scenarios influences emotions, intentions, decisions and moral judgments \citep{byrne2000counterfactual,epstude2008functional,buchsbaum2012power,van2015cognitive,byrne2016counterfactual,gerstenberg2024counterfactual}. As with all notions that occupy such a fundamental place in  human thought and affairs, counterfactuals warrant a rigorous, axiomatic mathematical foundation, to enable quantitative analyses and principled applications. Such is the goal of this paper.

Counterfactuals have been studied and formalised in a variety of ways across philosophy, logic, psychology, economics and computer science \citep{mandel2007psychology,heckman2007handbook,halpern2016actual}, with differing degrees of emphasis on stochasticity; in this work, we focus specifically on their stochastic aspects. 
Hence, we rely heavily on the axiomatisation of \emph{probability theory} due originally to \citet{kolmogorov1933foundations}, which has since become the widely accepted mathematisation of stochastics. One field where stochastic counterfactuals are prominently discussed is in the field of \emph{causality}, yet another cornerstone of human cognition, as well as of the sciences \citep{beebee2009oxford,illari2011causality,waldmann2017oxford}. 
%
%
Here, one is interested in studying the effects of \emph{interventions}, in contrast to passive observation of the world, as one does in (\say{pure}) probability theory. The relationship between causality and counterfactuals has been studied vigorously by philosophers and statisticians alike \citep{collins2004causation,pearl2000causal}: when considering the causal effect of an action, one compares, implicitly or explicitly, the ensuing events against those in an \say{imagined} world in which the original action is different. Many of the current mathematical theories of counterfactuals are founded upon a mathematical framework of causality, most saliently, those employing the so-called \emph{structural causal models} (SCMs) of \citet{pearl2009causality}, or the potential outcomes (POs) of \citet{rubin2005causal}. 

\begin{figure}[t]
    \centering
    \begin{tikzpicture}[roundnode/.style={circle, draw=black, very thick, minimum size=7mm}]
        \node[align=center] (association) at (0.5, 0) {Observation};
        \node (xassoc) at (-2.5,-1) {\(X\)};
        \node (yassoc) at (-1.2,-1) {\(Y\)};
        \draw[thick, ->](xassoc.east) to (yassoc.west);
        \node (uyassoc) at (-1.2,0) {\(U_Y\)};
        \node (uxassoc) at (-2.5,0) {\(U_X\)};
        \draw[thick, ->](uxassoc.south) to (xassoc.north);
        \draw[thick, ->](uyassoc.south) to (yassoc.north);
        \node[align=center] (intervention) at (0.5, 1.5) {Intervention};
        \node (xint) at (-2.5,1) {\(\text{do}(X)\)};
        \node (yint) at (-1.2,1) {\(Y\)};
        \draw[thick, ->](xint.east) to (yint.west);
        \node (uyint) at (-1.2,2) {\(U_Y\)};
        \node (uxint) at (-2.5,2) {\(U_X\)};
        \draw[thick, ->](uyint.south) to (yint.north);
        \draw[thick, dotted, ->](uxint.south) to (xint.north);
        \node[inner sep=0pt] (hammer) at (-3,1.2) {\includegraphics[scale=0.02]{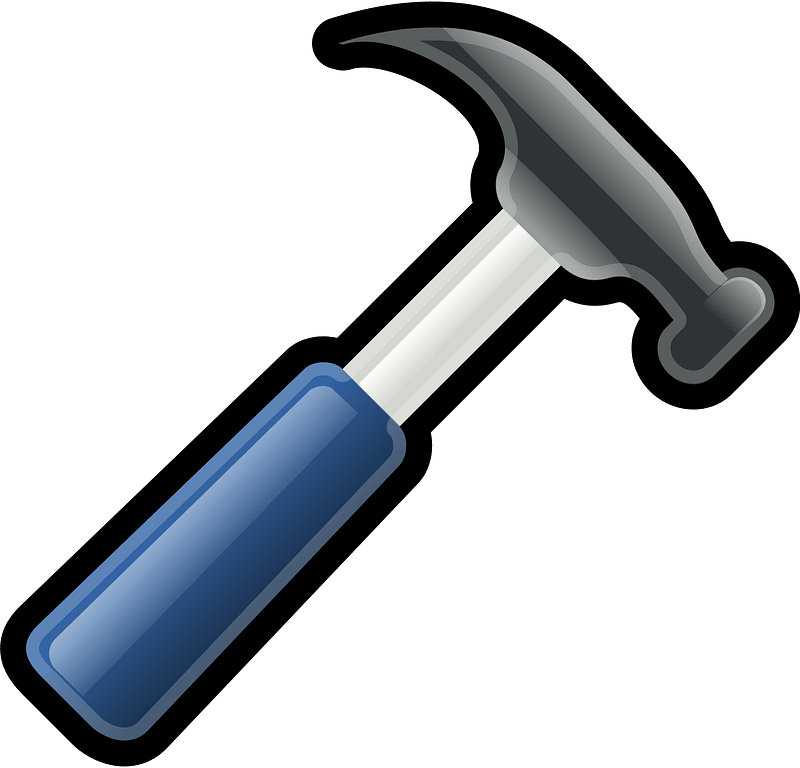}};
        \node (counterfactual) at (0.5, 3) {Counterfactual};
        \node (xcount) at (-2.5,3) {\(\text{do}(X)\)};
        \node (ycount) at (-1.2,3) {\(Y\)};
        \draw[thick, ->](xcount.east) to (ycount.west);
        \node (uycount) at (-1.2,4) {\(U_Y|X,Y\)};
        \node (uxcount) at (-2.5,4) {\(U_X\)};
        \draw[thick, dotted, ->](uxcount.south) to (xcount.north);
        \draw[thick, ->](uycount.south) to (ycount.north);
        \node[inner sep=0pt] (hammer2) at (-3,3.2) {\includegraphics[scale=0.02]{hammer.png}};
        \node[inner sep=0pt] (eyes) at (-1.5,4.3) {\includegraphics[scale=0.02]{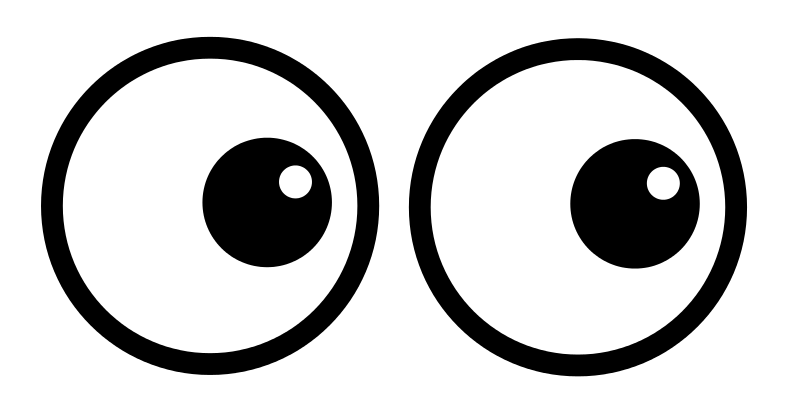}};
        \draw[thick, ->](association.north) to (intervention.south);
        \draw[thick, ->](intervention.north) to (counterfactual.south);
        \node[align=center] (probability) at (4,0.5) {Probability spaces\\\((\Omega,\sH,\PP)\)};
        \node[align=center] (causal) at (10,0.5) {Causal spaces\\\((\Omega,\sH,\PP,\KK)\)};
        \node[align=center] (counterfactual) at (4,2.5) {Counterfactual\\probability spaces\\\((\fact{\Omega}\times\cfact{\Omega},\fact{\sH}\otimes\cfact{\sH},\PP)\)};
        \node[align=center] (causalcounterfactual) at (10,2.5) {Counterfactual\\causal spaces\\\((\fact{\Omega}\times\cfact{\Omega},\fact{\sH}\otimes\cfact{\sH},\PP,\KK)\)};
        \draw[thick, ->](probability.east) to (causal.west);
        \draw[thick, ->](probability.north) to (counterfactual.south);
        \draw[thick, ->](counterfactual.east) to (causalcounterfactual.west);
        \draw[thick, ->](causal.north) to (causalcounterfactual.south);
    \end{tikzpicture}
    \caption{Left: Pearl's ladder of causation. Concepts in the upper rungs are strict generalisations of those in the lower rungs. SCMs are used to calculate the observational, interventional and counterfactual distributions in all of the rungs. Right: the view explored in this paper. Causal spaces and counterfactual probability spaces are each orthogonal extensions of probability spaces, and combining the two, we obtain counterfactual causal spaces. }
    \label{fig:ladder}
\end{figure}
The widely held view in the causality community, epitomised by Pearl's celebrated \say{ladder of causation} \citep{pearl2018book,bareinboim2022pearl}, is that interventions comprise a fundamental building block of counterfactuals (cf. \Cref{fig:ladder}, left). We explore a significant departure from these existing approaches that nest counterfactuals inside a causal model. 
While interventional counterfactuals are undoubtedly important, humans also often reason about counterfactual scenarios in which no intervention takes place in any of the worlds (e.g. \citealt{byrne2007precis}).
Therefore, the mathematisation of such non-interventional counterfactuals need not be based on a causal model that encodes interventional information. Contrary to this idea, all of the current approaches, even those variants that explicitly only treat non-interventional counterfactuals (for example, the extended SCMs of \citealt{lucas2015improved} and backtracking SCMs of \citealt{von2023backtracking}), base their mathematics on causal models designed for interventions. 

Further, existing frameworks of counterfactuals are often plagued by stringent assumptions inherent in their mathematisations (see \citet{park2023measure} for a detailed discussion discussion), such as acyclicity (see \Cref{ex:exam_cycle} for a case where acyclicity is not satisfied), discreteness (severely limiting the treatment of continuous-time stochastic processes), and that the endogenous variables do not causally affect the exogenous variables. The latter assumption is crucial in the abduction--action--prediction paradigm of counterfactuals in the SCM framework, but situations where it is not reasonable are ubiquitous: for example, a model of supply and demand in economics cannot be expected to have included all the variables that both affect and are affected by supply and demand. 

Finally, traditional approaches focus 
on the case in which as much is shared between the worlds as possible \emph{but} for a (typically small) chosen part of the system. This is the guiding principle behind the influential account of \say{similar worlds} by \citet{lewis1973counterfactuals}, as well as the SCM framework \citep{pearl2009causality,peters2017elements}. In the latter, the worlds share the same values for all of the noise variables, and the structural equations that are not intervened upon. Maximal similarity between worlds is not merely desired; it is stipulated by definition. 

In this paper, we propose an alternative perspective that views counterfactuals as an orthogonal concept to interventions (cf. \Cref{fig:ladder}, right). In particular, we argue that interventions are not a necessary ingredient for the formalisation of counterfactuals. Rather, we consider the incorporation of counterfactual outcomes and events as the essence of the study of counterfactuals, and this does not necessitate the introduction of wholly new mathematical objects. Accordingly, we define \emph{counterfactual probability spaces} as special cases of probability spaces by taking the product of two (or more) measurable spaces, each representing a \say{world}. In a similar spirit, \emph{counterfactual causal spaces} are defined as special cases of causal spaces, a recently proposed measure-theoretic axiomatisation of interventional causality \citep{park2023measure}, by requiring that the underlying measurable space be a product of world-specific component measurable spaces. We use \emph{counterfactual spaces} as an umbrella term to refer to both counterfactual probability and causal spaces. These definitions allow orthogonal mathematisations of interventions and counterfactuals (cf. \Cref{fig:ladder}, right). 
Further, because probability spaces and causal spaces are axiomatic formalisms that impose only minimal assumptions on the data-generating process, counterfactual spaces inherit the same level of generality. 



A key consequence of our mathematisation of stochastic counterfactuals is that it allows us to 
arbitrarily model how much information is shared between the worlds---in other words, how they are related to each other. 
Specifically, in counterfactual spaces, the similarity between worlds is encoded in the probability measure and the causal mechanism, the former representing the shared information in the observational state, and the latter that after interventions. Events in different worlds can be independent---meaning that there is no shared information---or synchronised---corresponding to maximal shared information---or anything in between. 

In this way, counterfactual spaces strictly generalise the existing frameworks, while being capable of incorporating a broader spectrum of counterfactuals. The type of counterfactuals typically considered in the usual SCM framework is conditioning in the factual world and intervening in the counterfactual world through the \say{abduction--action--prediction} procedure \citep{pearl2009causality}. Looking ahead to our running example of students attending a revision class and their exam results (\Cref{ex:exam}), the type of queries that can be answered using the above scheme in the usual SCMs is of the form,
\begin{quote}
    \say{Given that the student did not attend the class and failed, what is the probability that they would have passed if they had been \emph{forced} to attend the class?}
\end{quote}
On the other hand, backtracking SCMs \citep{lucas2015improved,von2023backtracking} are able to answer queries of the form,
\begin{quote}
    \say{Given that the student did not attend the class and failed, would they have passed if they had been \emph{observed} to attend the class?}
\end{quote}
Although above two queries appear similar, interventions and observations are fundamentally different, a distinction that underlies the entire concept/field of causality. 
POs, by contrast, consider counterfactual worlds that each corresponds to a hard intervention on the treatment variable, and answers queries on the joint distribution over the worlds, such as
\begin{quote}
    \say{What is the probability that the student passes if they attend the class and fails if they do not attend the class?}
\end{quote} 
Counterfactual spaces provides a unifying framework that allows one to answer all of the above queries and much more, by conditioning and intervening in either or both worlds in any combination, in any sequence, and with any amount of shared information between the worlds, before and after intervention. 

The paper is organised as follows. After introducing the necessary background on probability and causal spaces in \Cref{sec:preliminaries}, we define counterfactual probability spaces in \Cref{sec:counterfactual_probability_spaces}, and counterfactual causal spaces in \Cref{sec:counterfactual_causal_spaces}. In \Cref{sec:multiple}, we construct counterfactual spaces to incorporate more than two parallel worlds, and finally, in \Cref{sec:related_works}, we show that counterfactual spaces strictly generalise the SCM and PO frameworks, by explicitly constructing counterfactual spaces starting from arbitrary specifications of these frameworks.  

\section{Preliminaries \& notation}\label{sec:preliminaries}
In this section, we introduce the notation 
and 
recall the main concepts of probability and causal spaces that we will use in the manuscript. We require the former to define counterfactual probability spaces, and emphasize that the latter is only needed for counterfactual causal spaces.

\subsection{Probability theory}\label{subsec:probability_theory}
We first recall the axioms of probability theory. For a comprehensive introduction, see, for example, \citep{cinlar2011probability,durrett2019probability}. 
\begin{definition}\label{def:probability_space}
    A \emph{probability space} is a triple \((\Omega,\sH,\PP)\), where \(\Omega\) is a set of outcomes, \(\sH\) is a \(\sigma\)-algebra of events satisfying
    \begin{enumerate}[(i)]
        \item \(\Omega\in\sH\);
        \item \(A\in\sH\implies\Omega\setminus A\in\sH\);
        \item \(A_1,A_2,...\in\sH\implies\cup_{n=1}^\infty A_n\in\sH\);
    \end{enumerate}
    and \(\PP\) is a probability measure on \(\sH\), i.e. a function \(\PP:\sH\to[0,1]\) satisfying
    \begin{enumerate}[(i)]
        \item \(\PP(\Omega)=1\);
        \item \(\PP(\cup^\infty_{n=1}A_n)=\sum^\infty_{n=1}\PP(A_n)\) for every pairwise disjoint sequence \((A_n)\) in \(\sH\).
    \end{enumerate}
\end{definition}
We denote by \(\EE\) the expectation of a random variable with respect to the measure \(\PP\). For \(\omega\in\Omega\), the {Dirac measure} \(\delta_\omega:\sH\to\{0,1\}\) is defined such that \(\delta_\omega(A)=1\) if \(\omega\in A\) and \(0\) otherwise. Similarly, for any \(A\in\sH\), the indicator function \(\ind_A:\Omega\to\{0,1\}\) is defined such that \(\ind_A(\omega)=1\) if \(\omega\in A\) and \(0\) otherwise. 

For an event \(G\in\sH\), we denote by \(\PP_G\) the conditional probability given \(G\), defined, for each event \(A\in\sH\), by
\begin{equation*}
    \PP_G(A)=\begin{cases}\frac{\PP(G\cap A)}{\PP(G)}&\text{if }\PP(G)>0\\\text{undefined}&\text{otherwise}.\end{cases}
\end{equation*}
For any sub-\(\sigma\)-algebra \(\sG\) of \(\sH\), we denote by \(\PP_\sG\) the conditional probability given \(\sG\), defined, for each event \(A\in\sH\), as any \(\sG\)-measurable random variable \(\omega\mapsto\PP_\sG(\omega,A)\) such that, for any \(B\in\sG\), we have
\begin{equation*}
    \PP(A\cap B)=\EE\left[\ind_B(\cdot)\PP_\sG(\cdot,A)\right].
\end{equation*}
Throughout, we will use \(G\) and \(\sG\) for the event and \(\sigma\)-algebra to condition on. Note that \(\PP_G\) and \(\PP_\sG\) are different objects---for a fixed \(A\in\sH\), the former is a single positive number, whereas the latter is a \(\sG\)-measurable random variable, which can be shown to exist uniquely up to \(\PP\)-null events. For more details, see e.g. \citet[Chapter IV, Section 1]{cinlar2011probability}. 

Finally, we recall the definitions of (conditional) independence and almost sure equality of events. The former encodes the fact that no information is shared between events or \(\sigma\)-algebras, and the latter that maximal information is shared. These concepts will play an important role in our discussions of shared information between factual and counterfactual worlds. 
\begin{definition}\label{def:independence}
    Let us take a probability space \((\Omega,\sH,\PP)\), events \(A,B,G\in\sH\) and sub-\(\sigma\)-algebras \(\sF_1,\sF_2,\sG\subseteq\sH\).
    \begin{enumerate}[(i)]
        \item We say that \(A\) and \(B\) are \emph{independent}, and write \(A\indep{\PP} B\), if \(\PP(A\cap B)=\PP(A)\PP(B)\). 

        We say that \(\sF_1\) and \(\sF_2\) are \emph{independent}, and write \(\sF_1\indep{\PP}\sF_2\), if \(A\indep{\PP}B\) for all \(A\in\sF_1\) and all \(B\in\sF_2\). 
        \item We say that \(A\) and \(B\) are \emph{conditionally independent given \(G\)}, and write \(A\indep{\PP_G}B\), if \(\PP(G)>0\) and \(\PP_G(A\cap B)=\PP_G(A)\PP_G(B)\). 

        We say that \(\sF_1\) and \(\sF_2\) are \emph{conditionally independent given \(G\)}, and write \(\sF_1\indep{\PP_G}\sF_2\), if \(A\indep{\PP_G}B\) for all \(A\in\sF_1\) and all \(B\in\sF_2\). 
        \item We say that events \(A\) and \(B\) are \emph{conditionally independent given \(\sG\)}, and write \(A\indep{\PP_\sG}B\), if \(\PP_\sG(\omega,A\cap B)=\PP_\sG(\omega,A)\PP_\sG(\omega,B)\) for \(\PP\)-almost all \(\omega\in\Omega\). 

        We say that \(\sF_1\) and \(\sF_2\) are \emph{conditionally independent given \(\sG\)}, and write \(\sF_1\indep{\PP_\sG}\sF_2\), if \(A\indep{\PP_\sG}B\) for all \(A\in\sF_1\) and all \(B\in\sF_2\). 
    \end{enumerate}
\end{definition}
\begin{definition}\label{def:almost_sure_equality}
    Let us take a probability space \((\Omega,\sH,\PP)\) and events \(A,B,G\in\sH\). Let \(A\Delta B=(A\cup B)\setminus(A\cap B)\) be the symmetric difference of \(A\) and \(B\). We say that \(A\) and \(B\) are
    \begin{enumerate}[(i)]
        \item \emph{almost surely equal}, and write \(A\asequal{\PP}B\), if \(\PP(A\Delta B)=0\);
        \item \emph{almost surely equal given \(G\)}, and write \(A\asequal{\PP_G}B\), if \(\PP(G)>0\) and \(\PP_G(A\Delta B)=0\). 
    \end{enumerate}
\end{definition}
The analogue of \Cref{def:independence}(iii) for \Cref{def:almost_sure_equality}, i.e. \say{almost sure equality given \(\sG\)}, is redundant, since \(\PP_\sG(\omega,A\Delta B)=0\) for almost all \(\omega\in\Omega\) if and only if \(A\asequal{\PP}B\). 

\subsection{Causal spaces}\label{subsec:causal_spaces}
We also recall the definition of \emph{causal spaces} \citep{park2023measure}. Here, the key object is the \emph{transition probability kernel}. For measurable spaces \((E,\sE)\) and \((F,\sF)\), a mapping \(K:E\times\sF\to[0,1]\) is called a transition probability kernel from \((E,\sE)\) into \((F,\sF)\) if
\begin{itemize}
    \item the mapping \(x\mapsto K(x,B)\) is measurable for every set \(B\in\sF\), and
    \item the mapping \(B\mapsto K(x,B)\) is a probability measure on \((F,\sF)\) for every \(x\in E\).
\end{itemize}
Under extremely mild conditions, conditional measures are transition probability kernels \citep[p.150, Definition 2.4 \& p.151, Theorem 2.7]{cinlar2011probability}. 

We require that the measurable space be in product form. Let \(T\) be the index set of the product. Then taking, for each \(t\in T\), a set \(\Omega_t\) and a \(\sigma\)-algebra \(\sE_t\) on \(\Omega_t\), we have the product measurable space
\begin{equation*}
    (\Omega,\sH)=\otimes_{t\in T}(\Omega_t,\sE_t)=(\times_{t\in T}\Omega_t,\otimes_{t\in T}\sE_t).
\end{equation*}
Here, and in the rest of the paper, we use the notation \(\otimes\) for the product \(\sigma\)-algebra, and as a slight (and widespread) abuse of notation, we also use \(\otimes\) for the product of measurable spaces. 

For each \(S\subseteq T\), we denote by \(\sH_S\) the sub-\(\sigma\)-algebra of \(\sH\) generated by measurable rectangles \(\times_{t\in T}A_t\), where \(A_t\in\sE_t\) for all \(t\in T\), \(A_t=\Omega_t\) for all \(t\notin S\), and \(A_t\neq\Omega_t\) for only finitely many \(t\in S\). In particular, \(\sH_\emptyset=\{\emptyset,\Omega\}\) is the trivial sub-\(\sigma\)-algebra, and \(\sH_T=\sH\) is the full \(\sigma\)-algebra. Also, we denote by \(\Omega_S\) the subspace \(\times_{s\in S}\Omega_s\) of \(\Omega\), and for each \(\omega=(\omega_t)_{t\in T}\in\Omega\), we write \(\omega_S=(\omega_s)_{s\in S}\in\Omega_S\), where \(\omega_s\in\Omega_s\) for each \(s\in S\). 

A causal space is defined as follows. 
\begin{definition}[{\citep[Definition 2.2]{park2023measure}}]\label{def:causal_space}
    A \emph{causal space} is defined as the quadruple \((\Omega,\sH,\PP,\KK)\), where \((\Omega,\sH)\) is a measurable space with the above product structure, \(\PP\) is a probability measure on \((\Omega,\sH)\) and \(\KK=\{K_S:S\subseteq T\}\), called the \emph{causal mechanism}, is a collection of transition probability kernels \(K_S\) from \((\Omega,\sH_S)\) into \((\Omega,\sH)\), called the \emph{causal kernel on \(\sH_S\)}, that satisfy the following axioms:
    \begin{enumerate}[(i)]
        \item for all \(A\in\sH\) and \(\omega\in\Omega\), we have \(K_\emptyset(\omega,A)=\PP(A)\);
        \item for all \(\omega\in\Omega\), \(A\in\sH_S\) and \(B\in\sH\), we have \(K_S(\omega,A\cap B)=\ind_A(\omega)K_S(\omega,B)\).
    \end{enumerate}
\end{definition}
The probability measure \(\PP\) is the \say{observational measure}, and \(\KK\) encodes the causal information, along with the notion of \emph{interventions}, defined in the following.
\begin{definition}[{\citep[Definition 2.3]{park2023measure}}]\label{def:interventions_causal}
    Let us take a causal space \((\Omega,\sH,\PP,\KK)\), a subset \(U\subseteq T\) and a probability measure \(\QQ\) on \((\Omega,\sH_U)\). An \emph{intervention on \(\sH_U\) via \(\QQ\)} yields a new causal space \((\Omega,\sH,\iPP{U}{\QQ},\iKK{U}{\QQ})\), where the \emph{intervention measure} \(\iPP{U}{\QQ}\) is a probability measure on \((\Omega,\sH)\) defined, for \(A\in\sH\), by
    \begin{equation*}
        \iPP{U}{\QQ}(A)=\int\QQ(d\omega)K_U(\omega,A)
    \end{equation*}
    and \(\iKK{U}{\QQ}=\{\iK{U}{\QQ}{S}:S\subseteq T\}\) is the \emph{intervention causal mechanism} whose \emph{intervention causal kernels} are
    \begin{equation*}
        \iK{U}{\QQ}{S}(\omega_S,A)=\int\QQ(d\omega'_{U\backslash S})K_{S\cup U}((\omega_S,\omega'_{U\backslash S}),A).
    \end{equation*}
\end{definition}
Hence, causal kernels of the original causal space precisely encode what the new measure and new causal kernels will be after an intervention. We will denote by \(\iEE{U}{\QQ}\) the expectation with respect to \(\iPP{U}{\QQ}\). 

We also recall the definition of \emph{causal effects} in causal spaces, which will be crucial for an axiom of counterfactual causal spaces (\Cref{sec:counterfactual_causal_spaces}). 
\begin{definition}\label{def:causal_effect}[{\citet[Definition B.1]{park2023measure}}]
    Let us take a causal space \((\Omega,\sH,\PP,\KK)\) (cf. \Cref{def:causal_space}), an intervention set \(U\subseteq T\), an event \(A\in\sH\) and a sub-\(\sigma\)-algebra \(\sF\) of \(\sH\) (not necessarily of the form \(\sH_S\) for some \(S\subseteq T\)). 
    \begin{enumerate}[(i)]
        \item If \(K_S(\omega,A)=K_{S\setminus U}(\omega,A)\) for all \(S\in\sP(T)\) and all \(\omega\in\Omega\), then we say that \(\sH_U\) has \emph{no causal effect on \(A\)}. We say that \(\sH_U\) has \emph{no causal effect on \(\sF\)} if, for all \(A\in\sF\), \(\sH_U\) has no causal effect on \(A\).
        \item If there exists \(\omega\in\Omega\) such that \(K_U(\omega,A)\neq\PP(A)\), then we say that \(\sH_U\) has an \emph{active causal effect on \(A\)}. We say that \(\sH_U\) has an \emph{active causal effect on \(\sF\)} if \(\sH_U\) has an active causal effect on some \(A\in\sF\). 
        \item\label{dormantcausaleffect} Otherwise, we say that \(\sH_U\) has a \emph{dormant causal effect on \(A\)}. We say that \(\sH_U\) has a \emph{dormant causal effect on \(\sF\)} if \(\sH_U\) does not have an active causal effect on any event in \(\sF\) and there exists \(A\in\sF\) on which \(\sH_U\) has a dormant causal effect.
    \end{enumerate}
\end{definition}
It was shown in \citet[Remark B.2(a)]{park2023measure} that it is not possible for a \(\sigma\)-algebra \(\sH_U\) to have both no causal effect and an active causal effect on an event \(A\). But the definition of no causal effect is stronger than that of no active causal effect. No causal effect means that, not only does the measure of \(A\) remain the same as the observational measure \(\PP(A)\), but that intervening on \emph{any other \(\sigma\)-algebra} \(\sH_S\) is the same as intervening only on those components of \(S\) that do not belong to \(U\), i.e. on \(\sH_{S\setminus U}\). It is this stronger notion that we will need for counterfactual causal spaces. 

In the following simple toy example, we give an instantiation of a causal space and illustrate each of the above concepts of causal effect. 

\begin{example}\label{ex:dormant}
    Let us take three binary outcome sets \(\Omega_1=\Omega_2=\Omega_3=\{0,1\}\), so that the outcome set \(\Omega=\Omega_1\times\Omega_2\times\Omega_3\) has 8 elements. Let \(\PP\) be the uniform observational measure, and let us specify a subset of the causal kernels as in \Cref{tab:dormant}. The marginal observational measure on \(\Omega_3\) (first row of \Cref{tab:dormant}) is
    \begin{equation*}
        \PP(\omega_3=0)=\PP(\omega_3=1)=1/2.
    \end{equation*}
    Intervening on \(\sH_1\) with \(\omega_1=0\) or \(\omega_1=1\) keeps the measure on \(\sH_2\) and \(\sH_3\) uniform (second and third rows of \Cref{tab:dormant}):
    \begin{equation*}
        K_1(0,\{\omega_3=0\})=K_1(0,\{\omega_3=1\})=1/2=K_1(1,\{\omega_3=0\})=K_1(1,\{\omega_3=1\}),
    \end{equation*}
    and so \(\sH_1\) has no active causal effect on \(\sH_3\). Intervening with \(\omega_2=0\) has an active causal effect on \(\sH_3\), since
    \begin{equation*}
        K_2(0,\{\omega_3=0\})=1/4\neq1/2=\PP(\omega_3=0).
    \end{equation*}
    Finally, intervening with \(\omega_{1,2}=(0,0)\) has an active causal effect on \(\sH_3\), since
    \begin{equation*}
        K_{1,2}((0,0),\{\omega_3=0\})=1/8\neq1/2=\PP(\omega_3=0),
    \end{equation*}
    and in particular, as \(K_{1,2}((0,0),\{\omega_3=0\})=1/8\neq1/4=K_2(0,\{\omega_3=0\})\), we gather that \(\sH_1\) has a dormant causal effect on \(\{\omega_3=0\}\): the intervention \(\omega_1=0\) has no active causal effect by itself, but the joint intervention \(\omega_{1,2}=(0,0)\) is different to the intervention \(\omega_2=0\). For \(\sH_1\) to have no causal effect on \(\sH_3\), we would have required \(K_{1,2}((\omega_1,\omega_2),A)=K_2(\omega_2,A)\) for all \(\omega_1\in\Omega_1\), all \(\omega_2\in\Omega_2\) and all \(A\in\sH_3\). 
    
    \begin{table}[t]
        \centering
        \begin{tabular}{c|cccccccc}
            \multirow{2}{*}{Outcome} & \multicolumn{8}{c}{\(\omega=(\omega_1,\omega_2,\omega_3)\)} \\
            & \((0,0,0)\) & \((1,0,0)\) & \((0,1,0)\) & \((0,0,1)\) & \((1,1,0)\) & \((1,0,1)\) & \((0,1,1)\) & \((1,1,1)\) \\\hline
            \(\PP\) & \(1/8\) & \(1/8\) & \(1/8\) & \(1/8\) & \(1/8\) & \(1/8\) & \(1/8\) & \(1/8\) \\
            \(K_1(0,\cdot)\) & \(1/4\) & \(0\) & \(1/4\) & \(1/4\) & \(0\) & \(0\) & \(1/4\) & \(0\) \\
            \(K_1(1,\cdot)\) & \(0\) & \(1/4\) & \(0\) & \(0\) & \(1/4\) & \(1/4\) & \(0\) & \(1/4\) \\
            \(K_2(0,\cdot)\) & \(1/8\) & \(1/8\) & \(0\) & \(3/8\) & \(0\) & \(3/8\) & \(0\) & \(0\) \\
            \(K_{1,2}((0,0),\cdot)\) & \(1/8\) & \(0\) & \(0\) & \(7/8\) & \(0\) & \(0\) & \(0\) & \(0\) \\
        \end{tabular}
        \caption{In \Cref{ex:dormant}, \(\sH_1\) has a dormant causal effect on \(\sH_3\) and \(\sH_2\) has an active causal effect on \(\sH_3\).}
        \label{tab:dormant}
    \end{table}
\end{example}
As we see in \Cref{ex:dormant}, the concepts of no causal effect and dormant causal effect are dependent on what other variables (or components of the measurable space) are included in the causal space. On the other hand, active causal effect is model-invariant, as long as the \(\sigma\)-algebra that we want to intervene on and the event on which we are interested are included. In other words, the former are not invariant to \emph{marginalisation} \citep{park2025fine}, whereas the latter is. In \Cref{ex:dormant}, if we marginalised out the component \(\Omega_2\), then not only would \(\sH_1\) not have an active causal effect on \(\sH_3\), but it would now have no causal effect on \(\sH_3\). 

Finally, we recall the definition of (active) conditional causal effects \citep{park2025fine}. This definition will have particular relevance in prototypical counterfactual queries of the form \say{given an observation in the factual world, what causal effect would an intervention in the counterfactual world have had?} We write \(\iPP{U}{\delta_\omega}(\cdot)\) and \(K_U(\omega,\cdot)\) interchangeably---it is immediate from \Cref{def:interventions_causal} that they are the same measure. 
\begin{definition}\label{def:conditional_causal_effect}
    Let us take a causal space \((\Omega,\sH,\PP,\KK)\), an intervention set \(U\subseteq T\), events \(A,G\in\sH\) and a \(\sigma\)-algebra \(\sF\subseteq\sH\). We say that \(\sH_U\) has an \emph{active causal effect} on \(A\) conditioned on \(G\) if there exists some \(\omega\in\Omega\) such that \(\PP(G)>0\), \(K_U(\omega,G)>0\) and \(\iPP{U}{\delta_\omega}_G(A)\neq\PP_G(A)\). 
\end{definition}

\section{Counterfactual probability spaces}\label{sec:counterfactual_probability_spaces}
In this section, we formalise non-interventional counterfactual reasoning, by using probability spaces whose measurable spaces are products of factual and counterfactual measurable spaces. We call such probability spaces \emph{counterfactual probability spaces}. We first give the definition that accommodates two parallel worlds, and later generalise to multiple worlds in \Cref{sec:multiple}. 

We denote the set of factual outcomes by \(\fact{\Omega}\), and the set of counterfactual outcomes by \(\cfact{\Omega}\).\footnote{In this paper, we use superscripts to denote worlds, and subscripts to denote components of worlds---see, for example, the product space constructed for causal spaces before \Cref{def:causal_space}.\label{foot:scripts}} We equip \(\fact{\Omega}\) and \(\cfact{\Omega}\) with \(\sigma\)-algebras \(\fact{\sE}\) and \(\cfact{\sE}\) respectively. Then the entire measurable space \((\Omega,\sH)\) is obtained by taking the product of the factual and counterfactual measurable spaces as follows:
\begin{equation*}
    (\Omega,\sH)=(\fact{\Omega}\times\cfact{\Omega},\fact{\sE}\otimes\cfact{\sE}).
\end{equation*}
For any outcome \(\omega\in\Omega\), we denote by \(\fact{\omega}\) and \(\omega^\text{CF}\) the projections of \(\omega\) to \(\fact{\Omega}\) and \(\cfact{\Omega}\) respectively, so that \(\omega\) is decomposed as \(\omega=(\fact{\omega},\omega^\text{CF})\). 

We denote by \(\fact{\sH}\) the sub-\(\sigma\)-algebra of \(\sH\) consisting of measurable cylinders \(A\times\cfact{\Omega}\), with \(A\in\fact{\sE}\). Likewise, we denote by \(\cfact{\sH}\) the sub-\(\sigma\)-algebra of \(\sH\) consisting of measurable cylinders \(\fact{\Omega}\times B\) with \(B\in\cfact{\sE}\).\footnote{Note that we have \(\sH=\fact{\sE}\otimes\cfact{\sE}\), but \emph{not} \(\sH=\fact{\sH}\otimes\cfact{\sH}\). This means that \(\fact{\sH}\) and \(\cfact{\sH}\) are sub-\(\sigma\)-algebras of \(\sH\), so that events in \(\fact{\sH}\) or \(\cfact{\sH}\) also belong to \(\sH\), but \(\fact{\sE}\) and \(\cfact{\sE}\) are not sub-\(\sigma\)-algebras of \(\sH\). Of course, isomorphisms exist to this effect.} We refer to events in \(\fact{\sH}\) as \emph{factual events}, to those in \(\cfact{\sH}\) as \emph{counterfactual events}, and to those that belong to neither as \emph{cross-world events}. Only the trivial events \(\emptyset\) and \(\Omega\) belong to both \(\fact{\sH}\) and \(\cfact{\sH}\). We also refer to sub-\(\sigma\)-algebras of \(\fact{\sH}\) as \emph{factual \(\sigma\)-algebras}, to sub-\(\sigma\)-algebras of \(\cfact{\sH}\) as \emph{counterfactual \(\sigma\)-algebras}, and to those that are neither as \emph{cross-world \(\sigma\)-algebras.}

We are ready to define counterfactual probability spaces. 
\begin{definition}\label{def:counterfactual_probability_space}
	A \emph{counterfactual probability space} is defined as the triple \((\Omega,\sH,\PP)\), where \((\Omega,\sH)\) is a measurable space with the above product structure and \(\PP\) is a probability measure on \((\Omega,\sH)\).
\end{definition}
Mathematically speaking, counterfactual probability spaces are simply probability spaces, with just the additional requirement that the measurable space be a product of the factual and counterfactual measurable spaces. There is no mathematical asymmetry between the factual and counterfactual measurable spaces---the nomenclature is for convenience. 

The interpretation is as follows. The marginals of the measure \(\PP\) on factual and counterfactual events are simply the probabilities that they occur in their corresponding worlds. The measure on the cross-world events is more interesting, as it tells us how much information is shared between the two worlds. For each pair of events \(A\in\fact{\sH}\) and \(B\in\cfact{\sH}\), if the measure is such that \(A\) and \(B\) are independent (see \Cref{def:independence}), then there is no shared information between \(A\) and \(B\). At the other extreme, if \(A\) and \(B\) are almost surely equal (see \Cref{def:almost_sure_equality}), then information share is maximal. At the level of worlds:
\begin{enumerate}
    \item if \(\fact{\sH}\indep{\PP}\cfact{\sH}\), then there is no shared information between the factual and counterfactual worlds;
    \item if \(\fact{\sH}\asequal{\PP}\cfact{\sH}\), i.e. for every \(A\in\fact{\sH}\), there exists \(B\in\cfact{\sH}\) such that \(A\asequal{\PP}B\) and vice versa, then the information share is maximal. In other words, conditioning on one world fully determines the other world. 
\end{enumerate}
Intuitively, the closer the counterfactual world is to the factual one---e.g. differing only by a local modification or a short time horizon---the more shared information we expect between them, whereas more distant counterfactuals tend to yield weaker cross-world dependence. 

We now give some examples of counterfactual probability spaces. 
\begin{example}\label{ex:coin}
    Let us take \(\fact{\Omega}=\cfact{\Omega}=\{H,T\}\), with
    \begin{equation*}
        \PP(\{(H,H)\})=\PP(\{(H,T)\})=\PP(\{(T,H)\})=\PP(\{(T,T)\})=0.25.
    \end{equation*}
    In this example, one unbiased coin is being flipped, once in the factual world and once in the counterfactual world. The events in the two worlds are independent under the measure \(\PP\), i.e. \(\fact{\sH}\indep{\PP}\cfact{\sH}\), and there is no information shared between the two worlds. We could alternatively specify \(\PP\) such that
    \begin{equation*}
        \PP(\{(H,H)\})=\PP(\{(T,T)\})=0.5,\qquad\PP(\{(H,T)\})=\PP(\{(T,H)\})=0.
    \end{equation*}
    With this measure, no other randomness enters the counterfactual world than those already present in the factual world. 
\end{example}
\begin{example}\label{ex:disease}
    Let us take \(\fact{\Omega}=\cfact{\Omega}=\{S,D\}\), with
    \begin{equation*}
        \begin{matrix}
            &&\multicolumn{2}{c}{\text{Counterfactual}} \\
            &\PP & S & D \\
            \multirow{2}{*}{\text{Factual}} & S & 0.89 & 0.01 \\
            & D & 0.01 & 0.09
        \end{matrix}
    \end{equation*}
    Patients with a disease have a 90\% chance of surviving (\(S\)) and a 10\% chance of dying (\(D\)). If a patient was observed to survive in the factual world, the probability of this person surviving in the counterfactual world is \(\frac{0.89}{0.89+0.01}\approx0.99\). Likewise, if a patient was observed to die in the factual world, then they also die in the counterfactual world with probability \(\frac{0.09}{0.09+0.01}=0.9\).
    
    Here, the shared information is induced by a query about a randomly chosen patient with the same underlying health conditions across both worlds, but other sources of randomness that also influence the survival are not shared between the worlds. 
\end{example}
\begin{example}\label{ex:exam}
    Suppose that we want to model the probability of a student attending a revision class and passing a subsequent exam, in two parallel worlds:
    \begin{alignat*}{3}
        \fact{\Omega}_\class&=\cfact{\Omega}_{\class}=\{Y,N\},\qquad&\fact{\Omega}_{\exam}&=\cfact{\Omega}_{\exam}=\{P,F\}\\
        \fact{\Omega}&=\fact{\Omega}_{\class}\times\fact{\Omega}_{\exam},&\cfact{\Omega}&=\cfact{\Omega}_{\class}\times\cfact{\Omega}_{\exam},\qquad\Omega=\fact{\Omega}\times\cfact{\Omega},
    \end{alignat*}
    where \(Y\), \(N\), \(P\) and \(F\) respectively stand for outcomes \say{Yes}, \say{No}, \say{Pass} and \say{Fail}. The full measure \(\PP\) is given in \Cref{tab:exam_observation}.
    \begin{table}[t]
        \begin{center}
            \begin{tabular}{cccccc|c}
                &&\multicolumn{5}{c}{Counterfactual}\\
                &\(\PP\)&\((Y,P)\)&\((Y,F)\)&\((N,P)\)&\((N,F)\) & Sum \\
                \multirow{5}{*}{Factual } &\((Y,P)\)&0.32&0.04&0.06&0.01&0.43\\
                &\((Y,F)\)&0.04&0.12&0.01&0.04&0.21\\
                &\((N,P)\)&0.06&0.01&0.1&0.02&0.19\\
                &\((N,F)\)&0.01&0.04&0.02&0.1&0.17\\\cline{2-7}
                & \text{Sum} & 0.43 & 0.21 & 0.19 & 0.17 & 1\\
            \end{tabular}
        \end{center}
        \caption{The measure across the factual and counterfactual worlds on a student attending the class and passing the exam. }
        \label{tab:exam_observation}
    \end{table}
    Using this measure, we can answer \say{backtracking counterfactual} queries \citep{von2023backtracking}, for example: 
    \begin{enumerate}[(a)]
        \item \say{Given that a student passed the exam after attending the revision class, what is the probability that the same student passes the same exam had they sat it again?} To answer this question, we condition on the first row, and calculate the sum of the first and the third columns, to obtain \(\frac{0.32+0.06}{0.43}\approx0.88\), which is higher than the marginal probability of a student passing the exam (\(0.43+0.19=0.62\)). 
        \item \say{Given that a student attended the class, what is the probability that the same student will attend the class if we turned back time?} For this, we would condition on the first two rows, and calculate the sum of the first two columns: \(\frac{0.32+0.04+0.04+0.12}{0.43+0.21}=0.8125\). Again, this is higher than the marginal probability of a student attending a class (\(0.43+0.21=0.64\)). 
        \item \say{Given that a student failed the exam after not attending the class, would the same student have passed the same exam if they were observed to attend the class instead?} To answer this question, we condition on the last row and the first two columns, and look at the first column: \(\frac{0.01}{0.01+0.04}=0.2\). Note that this is still much lower than the marginal probability 0.62 of passing, which makes sense because the ability of the student and the difficulty of the exam remain the same in the counterfactual world. However, it is higher than the probability of the student passing after simply conditioning on the student not attending the class and failing in the factual world, which is obtained by conditioning on the last row and summing the first and third columns: \(\frac{0.01+0.02}{0.17}\approx0.176\). 
        
        Note that this is different to asking \say{if they had been forced to attend the class?}---an observation is different to an intervention. Consequently, the above discussion tells us nothing about the causal relationship between attending the class and passing the exam. For causality, we need (as we always do) the notion of interventions, which is not treated in counterfactual probability spaces. We will consider interventions in counterfactual causal spaces, in \Cref{sec:counterfactual_causal_spaces}, and revisit this example. 
    \end{enumerate}
\end{example}

Of course, these notions can be extended in a straightforward manner to conditional statements. 
We now give an example of a case in which we have conditional synchronisation of factual and counterfactual events. 
\begin{example}\label{ex:star}
    Suppose that we model the observation of a particular star on a specific night. We take
    \begin{alignat*}{3}
        \fact{\Omega}_\sky&=\cfact{\Omega}_\sky=\{C,O\},\qquad&\fact{\Omega}_\starr&=\cfact{\Omega}_\starr=\{Y,N\},\\
        \fact{\Omega}&=\fact{\Omega}_\sky\times\fact{\Omega}_\starr,&\cfact{\Omega}&=\cfact{\Omega}_\sky\times\cfact{\Omega}_\starr,\qquad\Omega=\fact{\Omega}\times\cfact{\Omega}
    \end{alignat*}
    where \(C\), \(O\), \(Y\) and \(N\) respectively stand for the outcomes \say{Clear}, \say{Overcast}, \say{Yes} and \say{No}. The full measure \(\PP\) is given in \Cref{tab:star}. 
    \begin{table}[t]
        \begin{center}
            \begin{tabular}{cccccc|c}
                &&\multicolumn{5}{c}{Counterfactual}\\
                &\(\PP\)&\((C,Y)\)&\((C,N)\)&\((O,Y)\)&\((O,N)\)&Sum\\
                \multirow{5}{*}{Factual } &\((C,Y)\)&0.2&0&0.05&0.15&0.4\\
                &\((C,N)\)&0&0.05&0&0.05&0.1\\
                &\((O,Y)\)&0.05&0&0.01&0.04&0.1\\
                &\((O,N)\)&0.15&0.05&0.04&0.16&0.4\\\cline{2-7}
                &\text{Sum}&0.4&0.1&0.1&0.4&1\\
            \end{tabular}
        \end{center}
        \caption{The measure across the factual and counterfactual worlds on the sky and the star being observed. }
        \label{tab:star}
    \end{table}
    We note the following:
    \begin{itemize}
        \item The sky is equally likely to be clear or overcast, and the sky in the factual world is independent from the sky in the counterfactual world.
        \item The telescope used to observe the star is shared between the worlds, and has a \(1/5\) chance of being faulty, but it is marginalised out of the model. If the sky is clear, then the star will be observed with a working telescope without fail, but will not be observed with a faulty telescope. If the sky is overcast, the star will be observed with probability 1/4, and with a faulty telescope, it will not be observed. 
        \item We can see that the events \(\{\fact{\omega}_\starr=Y\}\) and \(\{\cfact{\omega}_\starr=Y\}\) are not almost surely equal, since summing up the first and third rows of the last column gives us \(\PP(\fact{\omega}_\starr=Y,\cfact{\omega}_\starr=N)=0.15+0.04=0.19>0\). Let us define the event in which the sky is clear in both worlds: \(G=\{\fact{\omega}_\sky=C,\cfact{\omega}_\sky=C\}\). Then conditioned on \(G\) (i.e. looking at the upper-left block of \Cref{tab:star}), the events \(\{\fact{\omega}_\starr=Y\}\) and \(\{\cfact{\omega}_\starr=Y\}\) are almost surely equal. Since we only have one binary variable in each world under conditioning on \(G\), this means that \(\fact{\sH}\asequal{\PP_G}\cfact{\sH}\), i.e. the factual and counterfactual worlds are synchronised given that the sky is clear in both worlds. This mathematically encodes that, on a clear night, the only random factor that determines the observation of the star is the telescope, which is shared across the worlds. 
    \end{itemize}
\end{example}
In the above examples, it was explicitly stated what was common in the two worlds (nothing in \Cref{ex:coin}, the patient in \Cref{ex:disease}, the student and the exam in \Cref{ex:exam} and the telescope in \Cref{ex:star}), but this was purely for the clarity of explanation. Mathematically, the measure \(\PP\) encodes the shared information, and the actual entity that is shared need not be (and mathematically is not) made explicit. Further, the formalism is agnostic to \emph{time}. On the one hand, one can interpret the counterfactual world as rolling back time and running events again (as we did in backtracking queries in \Cref{ex:exam}). On the other hand, one could also interpret both worlds as taking place in the future, starting from a common time point. Then, the further away the two worlds are from this common starting point, the less information they share. 

The marginal measure on \(\fact{\sH}\) and \(\cfact{\sH}\) were identical in all the examples above. This, in general, need not be the case; in fact, \Cref{def:counterfactual_probability_space} even allows the measurable spaces \((\fact{\Omega},\fact{\sE})\) and \((\cfact{\Omega},\cfact{\sE})\) of the two worlds to be different. However, the special case of the two worlds being symmetric is of interest.
\begin{definition}\label{def:counterfactual_probability_space_symmetric}
    Let \((\Omega,\sH,\PP)\) be a counterfactual probability space. We say that \((\Omega,\sH,\PP)\) is symmetric if
    \begin{enumerate}[(a)]
        \item the two measurable spaces are the same, i.e. \((\fact{\Omega},\fact{\sE})=(\cfact{\Omega},\cfact{\sE})\);
        \item for any \(A,B\in\fact{\sE}=\cfact{\sE}\), we have \(\PP(A\times B)=\PP(B\times A)\). 
    \end{enumerate}
\end{definition}
The counterfactual probability spaces in \Cref{ex:coin,ex:disease,ex:exam} are easily seen to be symmetric. Modelling worlds to be symmetric makes sense when the two worlds have the same information (whether or not it is shared). Let us return to \Cref{ex:disease} but impose a different measure \(\PP\) that makes the worlds asymmetric.
\begin{example}\label{ex:disease_asymmetric}
    Suppose that in the counterfactual world, the healthcare system has collapsed, and patients with the disease are more likely to die. Accordingly, the measure is now:
    \begin{equation*}
        \begin{matrix}
            & & \multicolumn{2}{c}{\text{Counterfactual}} \\
            & \PP & S & D \\
            \multirow{2}{*}{\text{Factual}} & S & 0.6 & 0.3 \\
            & D & 0.001 & 0.099.
        \end{matrix}
    \end{equation*}
    The marginal measure in the factual world, where the healthcare system is intact, remains the same (\(90\%\) chance of survival and \(10\%\) chance of death). However, in the counterfactual world, the marginal measure of survival is only \(60.1\%\), and the marginal measure of death \(39.9\%\). The patient with the same underlying health conditions is still interpreted to be shared across the worlds, meaning that if they were observed to survive in the factual world, they are still more likely than not to survive in the counterfactual world (with probability \(\frac{0.6}{0.6+0.3}=\frac{2}{3}\)), despite the collapsed healthcare system. In particular, it is higher than the marginal probability 0.601 of survival in the counterfactual world. 
\end{example}

\section{Counterfactual causal spaces}\label{sec:counterfactual_causal_spaces}
In this section, we define \emph{counterfactual causal spaces} as special cases of causal spaces, whose measurable spaces are products of factual and counterfactual measurable spaces. This is analogous to how we obtained counterfactual probability spaces in \Cref{sec:counterfactual_probability_spaces}, as probability spaces whose measurable spaces are products of factual and counterfactual components. 
In addition, we also impose an extra axiom that there be no cross-world causal effect. In \Cref{subsec:counterfactual_causal_space}, we formally introduce counterfactual causal spaces, and interventions therein. In \Cref{subsec:synchronisation_independence_causal}, we again discuss the two extremes of shared information, namely, independence and synchronisation of worlds, in the context of counterfactual causal spaces. 

\subsection{Formal definition}\label{subsec:counterfactual_causal_space}
We first construct the underlying measurable space. Similarly as in causal spaces (\Cref{def:causal_space}), we require that the factual and counterfactual measurable spaces be in product form. We denote by \(\fact{T}\) and \(\cfact{T}\) the factual and counterfactual index sets, and write \(T=\fact{T}\cup\cfact{T}\). For each \(t\in\fact{T}\), we take a measurable space \((\fact{\Omega}_t,\fact{\sE}_t)\), and for each \(t\in\cfact{T}\), we take \((\cfact{\Omega}_t,\cfact{\sE}_t)\). Then, we define the sets of factual and counterfactual outcomes respectively as \(\fact{\Omega}=\times_{t\in \fact{T}}\fact{\Omega}_t\) and \(\cfact{\Omega}=\times_{t\in \cfact{T}}\cfact{\Omega}_t\). Also, denote by \(\fact{\sE}=\otimes_{t\in \fact{T}}\fact{\sE}_t\) and \(\cfact{\sE}=\otimes_{t\in \cfact{T}}\cfact{\sE}_t\) the corresponding \(\sigma\)-algebras. Then the entire measurable space \((\Omega,\sH)\) is obtained by taking the product of the factual and counterfactual measurable spaces as follows:
\begin{equation*}
    (\Omega,\sH)=(\fact{\Omega}\times\cfact{\Omega},\fact{\sE}\otimes\cfact{\sE}).
\end{equation*}
Similarly as in \Cref{sec:counterfactual_probability_spaces}, for any outcome \(\omega\in\Omega\), we denote by \(\fact{\omega}\) and \(\cfact{\omega}\) its projections to \(\fact{\Omega}\) and \(\cfact{\Omega}\) respectively, so that \(\omega\) is decomposed as \(\omega=(\fact{\omega},\cfact{\omega})\). Further, similarly as in \Cref{subsec:causal_spaces}, for any \(S\subseteq T\), we denote by \(\Omega_S\) the subspace \(\times_{s\in S}\Omega_s\) of \(\Omega\). We also write \(\omega_S=(\omega_s)_{s\in S}\), and if \(S\subseteq\fact{T}\) (respectively \(S\subseteq\cfact{T}\)), we also write \(\fact{\omega}_S\) (respectively \(\cfact{\omega}_S\)). 

For any \(S\subseteq T\), we denote by \(\sH_S\) the sub-\(\sigma\)-algebra of \(\sH\) generated by measurable rectangles \((\times_{t\in \fact{T}}A_t)\times(\times_{t\in \cfact{T}}B_t)\), where \(A_t\in\fact{\sE}_t\) and \(B_t\in\cfact{\sE}_t\) differ from \(\fact{\Omega}_t\) and \(\cfact{\Omega}_t\) only for finitely many \(t\) such that \(t\in S\). As a shorthand, we write \(\fact{\sH}=\sH_{\fact{T}}\) and \(\cfact{\sH}=\sH_{\cfact{T}}\). Just as in \Cref{sec:counterfactual_probability_spaces}, we refer to events in \(\fact{\sH}\) as \emph{factual events}, to those in \(\cfact{\sH}\) as \emph{counterfactual events}, and to those that belong to neither as \emph{cross-world events}. We also refer to sub-\(\sigma\)-algebras of \(\fact{\sH}\) as \emph{factual \(\sigma\)-algebras}, to sub-\(\sigma\)-algebras of \(\cfact{\sH}\) as \emph{counterfactual \(\sigma\)-algebras}, and to those that are neither as \emph{cross-world \(\sigma\)-algebras.}

We are finally ready to define counterfactual causal spaces. 
\begin{definition}\label{def:counterfactual_causal_space}
    A \emph{counterfactual causal space} is a quadruple \((\Omega,\sH,\PP,\KK)\), where \((\Omega,\sH)\) is a measurable space with the above product structure, \(\PP\) is a probability measure on \((\Omega,\sH)\) 
    and \(\KK=\{K_S:S\subseteq T\}\), called the \emph{causal mechanism}, is a collection of transition probability kernels \(K_S\) from \((\Omega,\sH_S)\) into \((\Omega,\sH)\), called the \emph{causal kernel on }\(\sH_S\), satisfying the following axioms:
	\begin{enumerate}[(i)]
		\item\label{item:trivial_intervention} for all \(\omega\in\Omega\) and \(A\in\sH\), we have
        \begin{equation*}
            K_\emptyset(\omega,A)=\PP(A);
        \end{equation*}
        \item\label{item:no_cross_world_effect} for all \(\omega\in\Omega\), all \(S\in\sP(T)\) and all \(A\in\fact{\sH}\), we have
        \begin{equation*}
            K_S(\omega,A)=K_{S\cap \fact{T}}(\omega,A),
        \end{equation*}
        and likewise, for all \(\omega\in\Omega\), all \(S\in\sP(T)\) and all \(B\in\cfact{\sH}\), we have
        \begin{equation*}
            K_S(\omega,B)=K_{S\cap \cfact{T}}(\omega,B);
        \end{equation*}
        \item\label{item:interventional_determinism} for all \(\omega\in\Omega\), \(A\in\sH^S\) and \(B\in\sH\), we have
        \begin{equation*}
            K_S(\omega,A\cap B)=\ind_A(\omega)K_S(\omega,B)=\delta_\omega(A)K_S(\omega,B);
        \end{equation*}
        in particular, for \(A\in\sH^S\), we have
        \begin{equation*}
            K_S(\omega,A)=\ind_A(\omega)K_S(\omega,\Omega)=\ind_A(\omega).
        \end{equation*}
    \end{enumerate}
\end{definition}
We will give intuitions on this definition and the axioms immediately after defining \emph{interventions} in counterfactual causal spaces:
\begin{definition}\label{def:intervention_counterfactual}
    Let us take a counterfactual causal space \((\Omega,\sH,\PP,\KK)\) (\Cref{def:counterfactual_causal_space}), an intervention set \(U\subseteq T\) and a probability measure \(\QQ\) on \((\Omega,\sH_U)\). 
    An \emph{intervention on \(\sH_U\) via \(\QQ\)} yields a new counterfactual causal space \((\Omega,\sH,\iPP{U}{\QQ},\iKK{U}{\QQ})\), where the \emph{intervention measure} \(\iPP{U}{\QQ}\) is a probability measure on \((\Omega,\sH)\) defined, for \(A\in\sH\), by
    \begin{equation*}
        \iPP{U}{\QQ}(A)=\int\QQ(d\omega)K_U(\omega,A)
    \end{equation*}
    and \(\iKK{U}{\QQ}=\{\iK{U}{\QQ}{S}:S\subseteq T\}\) is the \emph{intervention causal mechanism} whose \emph{intervention causal kernels} are
    \begin{equation*}
        \iK{U}{\QQ}{S}(\omega_S,A)=\int\QQ(d\omega'_{U\setminus S})K_{S\cup U}((\omega_S,\omega'_{U\setminus S}),A).
    \end{equation*}
\end{definition}
In the following remark, we give intuitions about the axioms of causal kernels in counterfactual causal spaces, given in \Cref{def:counterfactual_causal_space}. 
\begin{remark}\label{rem:axioms}
    Let \((\Omega,\sH,\PP,\KK)\) be a counterfactual causal space. 
    \begin{itemize}
        \item Axiom~\ref{item:trivial_intervention} tells us that if we do not intervene on anything, then the measure will stay the same as the initial measure \(\PP\).
        \item Axiom~\ref{item:no_cross_world_effect} tells us that there can be no cross-world causal effect (in the sense of \Cref{def:causal_effect}), i.e. factual \(\sigma\)-algebras have no causal effect on the counterfactual events, and vice versa. 
        \item Axiom~\ref{item:interventional_determinism} tells us that, after an intervention, the restriction of the resulting measure on the \(\sigma\)-algebra on which we intervened should coincide with the measure with which we intervened. 
    \end{itemize}
\end{remark}
Axioms~\ref{item:trivial_intervention} and \ref{item:interventional_determinism} are precisely the same as those of causal spaces (\Cref{def:causal_space}). Hence, just as counterfactual probability spaces were special cases of probability spaces, counterfactual causal spaces are special cases of causal spaces, but with an extra axiom which is not present in causal spaces. Moreover, of course, counterfactual causal spaces can be viewed as counterfactual probability spaces, by ignoring the causal mechanism. 

We must check that the counterfactual causal space obtained after an intervention is indeed a counterfactual causal space, i.e. \((\Omega,\sH,\iPP{U}{\QQ},\iKK{U}{\QQ})\) satisfies the axioms of \Cref{def:counterfactual_causal_space}. The proof of this statement is given as a special case of \Cref{thm:Nway_interventions}, where we prove the analogous result for \(N\)-way counterfactual causal spaces. 
\begin{savetheorem}{thm:intervention}
    The intervention causal mechanism \(\iKK{U}{\QQ}\) satisfies the axioms of \Cref{def:counterfactual_causal_space}. 
\end{savetheorem}

Let us make a few further remarks on counterfactual causal spaces and interventions. 
\begin{remark}\label{rem:counterfactual_causal_spaces}
    \begin{enumerate}[(i)]
        \item Again, there is no mathematical asymmetry between factual and counterfactual worlds---the nomenclature is for convenience and intuition. 
        \item Axiom~\ref{item:no_cross_world_effect} does not tell us about causal effects on cross-world events. Indeed, this is precisely how shared information after an intervention is encoded, which can be different to the information shared between the worlds before the intervention.  
        \item Marginalising a counterfactual causal space yields another counterfactual causal space, as long as an entire world is not marginalised out. This is because it is immediate that, if a \(\sigma\)-algebra \(\sH_U\) has no causal effect on an event \(A\) in the larger counterfactual causal space, it will have no causal effect on \(A\) in the smaller space. Hence, the no cross-world causal effect axiom (\Cref{def:counterfactual_causal_space}(ii)) is also preserved, and so the result of a marginalisation procedure is another counterfactual causal space.  
    \end{enumerate}
\end{remark}
Finally, we define \emph{symmetric} counterfactual causal spaces, analogously to symmetric counterfactual probability spaces (\Cref{def:counterfactual_probability_space_symmetric}). Here, not only do we require the probability measure \(\PP\) to be symmetric, but also all of the causal kernels. 
\begin{definition}\label{def:counterfactual_causal_space_symmetric}
    Let us take a counterfactual causal space \((\Omega,\sH,\PP,\KK)\), as defined in \Cref{def:counterfactual_causal_space}. We say that \((\Omega,\sH,\PP,\KK)\) is symmetric if
    \begin{enumerate}[(a)]
        \item the two index sets are the same, i.e. \(\fact{T}=\cfact{T}\), and each of the measurable sets are the same, i.e. for all \(t\in\fact{T}=\cfact{T}\), we have \((\fact{\Omega}_t,\fact{\sE}_t)=(\cfact{\Omega}_t,\cfact{\sE}_t)\) (this implies that \((\fact{\Omega},\fact{\sE})=(\cfact{\Omega},\cfact{\sE})\));
        \item for any events \(A,B\in\fact{\sE}=\cfact{\sE}\), we have \(\PP(A\times B)=\PP(B\times A)\);
        \item for any events \(A,B\in\fact{\sE}=\cfact{\sE}\), any subsets \(S\subseteq\fact{T}\), \(S'\subseteq\cfact{T}\) and any outcome \(\omega=(\omega',\omega'')\in\Omega\), we have \(K_{S\cup S'}((\omega',\omega''),A\times B)=K_{S'\cup S}((\omega'',\omega'),B\times A)\). 
    \end{enumerate}
\end{definition}
Depending on the intervention we carry out, a symmetric counterfactual causal space does not necessarily remain symmetric after an intervention. For example, if an intervention is carried out in only one of the worlds, then the subsequent counterfactual causal space is clearly not symmetric in general. Moreover, a symmetric counterfactual causal space is symmetric after marginalisation if and only if the same components of the measurable space are marginalised out in each world. 

Let us give an example of a counterfactual causal space, endowing the counterfactual probability space in \Cref{ex:exam} with a causal mechanism.
\begin{example}\label{ex:exam_counterfactual_causal_space}
    We specify (a subset of) the causal mechanism on the measurable space. The causal kernel \(K_{\cfact{\class}}\) corresponding to the interventions of making a student attend (\(Y\)) or not attend (\(N\)) a class in the counterfactual world is specified in \Cref{tab:exam_class,tab:exam_noclass}. We can specify \(K_{\fact{\class}}\) to be symmetrical, i.e. the transposes of \Cref{tab:exam_class,tab:exam_noclass}. 
    
    \begin{table}[t]
        \centering
        \begin{tabular}{cccccc|c}
            &&\multicolumn{5}{c}{Counterfactual}\\
            &\(K_{\cfact{\class}}(Y,\cdot)\)&\((Y,P)\)&\((Y,F)\)&\((N,P)\)&\((N,F)\) & Sum \\
            \multirow{5}{*}{Factual }&\((Y,P)\)&0.39&0.04&0&0&0.43\\
            &\((Y,F)\)&0.05&0.16&0&0&0.21\\
            &\((N,P)\)&0.16&0.03&0&0&0.19\\
            &\((N,F)\)&0.04&0.13&0&0&0.17\\\cline{2-7}
            &\text{Sum}&0.64&0.36&0&0&1\\
        \end{tabular}
        \caption{The causal kernel for intervening on the student to attend the class in the counterfactual world.}
        \label{tab:exam_class}
    \end{table}
    \begin{table}
        \centering
        \begin{tabular}{cccccc|c}
            &&\multicolumn{5}{c}{Counterfactual}\\
            &\(K_{\cfact{\class}}(N,\cdot)\)&\((Y,P)\)&\((Y,F)\)&\((N,P)\)&\((N,F)\) & Sum \\
            \multirow{5}{*}{Factual } &\((Y,P)\)&0&0&0.37&0.06&0.43\\
            &\((Y,F)\)&0&0&0.05&0.16&0.21\\
            &\((N,P)\)&0&0&0.15&0.04&0.19\\
            &\((N,F)\)&0&0&0.03&0.14&0.17\\\cline{2-7}
            & \text{Sum}&0&0&0.6&0.4& 1\\
        \end{tabular}
        \caption{The causal kernel for intervening on the student \emph{not} to attend the class in the counterfactual world.}
        \label{tab:exam_noclass}
    \end{table}
    Note first that, in accordance with the interventional determinism axiom (\Cref{def:counterfactual_causal_space}(iii)), the event of a student not attending the class (resp. attending the class) in the counterfactual world after intervening on them to attend (resp. not attend) the class has measure zero. Note also that, by the no cross-world causal effect axiom (\Cref{def:counterfactual_causal_space}(ii)), the marginal measure on the factual events remains the same as the marginal observational measure (the \say{Sum} column). 
    
    The full specification of the causal mechanism \(\KK\) would involve many more causal kernels, such as those corresponding to intervening on the exam result in either the factual or the counterfactual world (e.g. \(K_{\cfact{\exam}}(P,\cdot)\); see \Cref{ex:exam_cycle,tab:exam_pass}), or any combination of the two variables across the two worlds (e.g. \(K_{\cfact{\class},\fact{\exam}}(\{Y,F\},\cdot)\), etc.). In particular, the causal kernels in \Cref{tab:exam_class,tab:exam_noclass} only show that there is no active causal effect across worlds. In order to satisfy the no cross-world causal effect axiom (\Cref{def:counterfactual_causal_space}(ii)), the kernels corresponding to intervening in both worlds must be constructed to satisfy this axiom, e.g. \(K_{\fact{\class},\cfact{\class}}(\{Y,Y\},A)=K_{\cfact{\class}}(Y,A)\) for all counterfactual events \(A\). 
    
    With these causal kernels in hand, we can read off the tables (the \say{Sum} row at the bottom) that the probability of passing after intervening to make the student attend the class in the counterfactual world is 0.64, which is slightly higher than the marginal observational probability of passing, \(0.43+0.19=0.62\). Similarly, the probability of passing after intervening to prevent the student from attending the class is 0.6, which is slightly lower the marginal observational probability of passing. According to \Cref{def:causal_effect}(ii), this means that \(\sH_{\cfact{\class}}\) has an active causal effect on the event that the student passes the exam in the counterfactual world. 

    Let us have a look at a few queries that we can answer in this counterfactual causal space. We place a particular emphasis on the question of \emph{conditional causal effect} (\Cref{def:conditional_causal_effect})---given an observation in the factual world, we ask whether an intervention in the counterfactual world would have had a causal effect.
    \begin{enumerate}[(a)]
        \item \say{Given that a student did not attend the revision class and failed the exam, what would have been their probability of passing had the student been forced to attend the revision class?} To answer this, we condition on the last row of \Cref{tab:exam_class}: \(\frac{0.04}{0.04+0.13}\approx0.24\). This is still lower than the marginal observational probability of a student passing (0.62), but higher than the probability that the same student would have passed the same exam had they been left to make their own choice about attending the revision class (\(\frac{0.01+0.02}{0.17}\approx0.176\), the last row and the first and third columns of \Cref{tab:exam_observation}). 

        According to \Cref{def:conditional_causal_effect}, the above calculations mean that \(\cfact{\omega}_\class=Y\) has an active causal effect on the event \say{student passes the exam in the counterfactual world} conditioned on the observation that the student did not attend the revision class and failed the exam in the factual world. 
        \item \say{Given that a student passed the exam, what would have been their probability of passing had the student been prevented from attending the class?} We condition on the first and third rows of \Cref{tab:exam_noclass}, and sum the third column: \(\frac{0.37+0.15}{0.43+0.19}\approx0.838\). This is still much higher than the marginal observational probability of 0.62, reflecting the fact that a student who was capable of passing in the factual world is likely to pass again even if they cannot go to the revision class. However, it is slightly lower than the observational probability of passing conditioned on the student passing the exam (\(\frac{0.32+0.06+0.06+0.1}{0.43+0.19}\approx0.87\)), meaning the student was left to make their own choice about attending the revision class in the counterfactual world. Lastly, if the student was observed to pass in the factual world and was forced to go to the revision class in the counterfactual world, then the probability of passing, calculated by conditioning on the first and third rows of \Cref{tab:exam_class} and summing the first column, would be \(\frac{0.39+0.16}{0.43+0.19}\approx0.89\)---slightly higher still. 

        According to \Cref{def:conditional_causal_effect}, the above definition tells us that both \(\cfact{\omega}_\class=Y\) and \(\cfact{\omega}_\class=N\) have active causal effects on the event \say{student passes in the counterfactual world} conditioned on the event \say{student passes in the factual world}. 
        \item After observing, for example, that a student attends the class and passes the exam in the factual world, instead of asking what would have happened if they been prevented from attending the class, we can also ask what would have happened if they were forced to attend the class in the counterfactual world. At first glance, it may appear that the probability of the student passing should be the same as if we had not intervened at all in the counterfactual world---after all, the student attends the class in both worlds. However, unlike the SCM framework, observing that a student attends the class in the factual world does not, in general, guarantee that the student will attend the class again in the counterfactual world, even if we do not explicitly intervene so that the student does not attend the class in the counterfactual world. In other words, intervening to make the student attend the class in the counterfactual world after observing that the student attended the class (and passed) in the factual world can still have a (conditional) causal effect on the exam result. 

        Indeed, conditioning on the first rows of \Cref{tab:exam_observation,tab:exam_class}, we can see that, in the first case, the probability of passing is \(\frac{0.32+0.06}{0.43}\approx0.88\), whereas in the latter case, the probability of passing is \(\frac{0.39}{0.43}\approx0.91\). So according to \Cref{def:conditional_causal_effect}, the outcome \(\cfact{\omega_{\class}}=Y\) has a causal effect on the event that the student passes in the counterfactual world, conditioned on the event that the student attends the class and passes in the factual world. 
        \item On the other hand, suppose that the student was observed to be missing at the revision class and failed in the factual world. We can condition on the last row of \Cref{tab:exam_observation} to see that, conditioned on this observation in the factual world, the probability that they will pass in the counterfactual world is \(\frac{0.01+0.02}{0.17}\approx0.18\). If we had further intervened to prevent the student from attending the class in the counterfactual world, the probability of passing can be read off \Cref{tab:exam_noclass}, by conditioning on the last row again: \(\frac{0.03}{0.17}\approx0.18\). So in this case, these two probabilities are the same. This means that, according to \Cref{def:conditional_causal_effect}, letting \(G\) be the event that the student does not attend the class and fails the exam in the factual world, and \(A\) the event that the student passes in the counterfactual world, the outcome \(\cfact{\omega_\class}=N\) has no active causal effect on \(A\), conditioned on \(G\). 

        However, if we instead intervened to make the student attend the class in the counterfactual world, then we can condition on the last row of \Cref{tab:exam_class} to see that the probability of \(A\) is \(\frac{0.04}{0.17}\approx0.24\). Hence, the outcome \(\cfact{\omega_\class}=Y\) does have an active causal effect on \(A\) conditioned on \(G\),. 
    \end{enumerate}
\end{example}
We remark that the definitions of (conditional) causal effects in \Cref{def:causal_effect,def:conditional_causal_effect} are given as binary statements, i.e. whether or not there is \emph{any} causal effect \emph{at all}. It is out of the scope of this paper to discuss the \emph{nature and strength} of a causal effect, but we can see in (b) above that the (conditional) causal effect of attending the class conditioned on the student passing in the factual world, while present, is very small. 

It should also be noted that there may be cross-world conditional causal effects in counterfactual causal spaces. This may be surprising at first, since, in counterfactual causal spaces, the causal mechanism is \emph{axiomatically required} not to have any cross-world causal effects (\Cref{def:counterfactual_causal_space}(ii)). But it is only natural that this is so, because an intervention in the factual world may create, destroy or change the nature and/or strength of the shared information between the worlds. Mathematically speaking, let \(U\subseteq\fact{T}\), so that \(\sH_U\) is a factual \(\sigma\)-algebra, and let \(A\in\cfact{\sH}\) be a counterfactual event. Then \(\sH_U\) cannot have any causal effect on \(A\), but if it has a causal effect on \(G\), then \(\sH_U\) does have a conditional causal effect on \(A\) given \(G\). As we can see in (b) above, the observational probability of passing the exam in counterfactual world does not change after an intervention on class in the factual world, but conditioning on the exam outcome, the same intervention does affect the exam result in the counterfactual world. 

Of course, if we intervene in one world and condition only in the other world, then there cannot be any cross-world conditional causal effects. The proof is in \Cref{sec:proofs}. 
\begin{saveprop}{prop:no_cross_world_conditional}
    Let us take a counterfactual causal space \((\Omega,\sH,\PP,\KK)\), a subset \(U\subseteq\fact{T}\) (so that \(\sH_U\) is a factual \(\sigma\)-algebra) and counterfactual events \(A,G\in\cfact{\sH}\). Then \(\sH_U\) has no causal effect on \(A\) conditioned on \(G\).
\end{saveprop}
Clearly, this result also holds vice versa---if the intervention takes place in the counterfactual world, the conditioning takes place in the factual world and we are interested in a factual event. 

\subsection{Synchronisation and independence of worlds}\label{subsec:synchronisation_independence_causal}
Counterfactual causal spaces can be viewed as counterfactual probability spaces by ignoring the causal mechanism, so the definitions of (conditional) independence of the factual and counterfactual worlds and their being (conditionally) synchronised carry over from \Cref{sec:counterfactual_probability_spaces}. We can further define their analogues with the causal kernels, which will represent the two extremes of information shared between the worlds \emph{after} an intervention. 

We first recall the notion of \emph{causal independence}, which is a direct interventional analogue of conditional independence. 
\begin{definition}\label{def:causal_independence}[{\citet[Definition 3.4]{buchholz2024products}}]
    Let us take a causal space \((\Omega,\sH,\PP,\KK)\) as in \Cref{def:causal_space}. Then for \(U\subseteq T\), two events \(A,B\in\sH\) are \emph{causally independent on \(\sH_U\)}, and write \(A\indep{K_U}B\), if, for all \(\omega\in\Omega\), 
    \[K_U(\omega,A\cap B)=K_U(\omega,A)K_U(\omega,B).\]
    We say that two sub-\(\sigma\)-algebras \(\sF_1\) and \(\sF_2\) are \emph{causally independent on \(\sH_U\)}, and write \(\sF_1\indep{K_U}\sF_2\), if \(A\indep{K_U}B\) for all \(A\in\sF_1\) and \(B\in\sF_2\). 
\end{definition}
Note that, unlike conditional independence, we require the above property to hold \emph{for all} \(\omega\in\Omega\) for causal independence, not just almost surely. This is because, during an intervention, it is possible to impose a measure on \(\sH_U\) that gives positive measure on events that previously had zero measure. 

We also define a causal analogue of almost surely equal events, determination and synchronisation of \(\sigma\)-algebras (c.f. \Cref{def:almost_sure_equality}).
\begin{definition}\label{def:causally_equal}
    Let us take a causal space \((\Omega,\sH,\PP,\KK)\), events \(A,B\in\sH\), sub-\(\sigma\)-algebras \(\sF_1,\sF_2\subseteq\sH\) and a subset \(U\subseteq T\). 
    \begin{enumerate}[(i)]
        \item We say that \(A\) and \(B\) are \emph{causally equal} on \(\sH_U\), and write \(A\asequal{K_U}B\), if \(K_U(\omega,A\Delta B)=0\) for all \(\omega\in\Omega\). 
        \item We say that \(\sF_1\) and \(\sF_2\) are \emph{causally synchronised on \(\sH_U\)}, and write \(\sF_1\asequal{K_U}\sF_2\), if, for each \(A\in\sF_1\), there exists \(B\in\sF_2\) such that \(A\asequal{K_U}B\) and vice versa.  
    \end{enumerate}
\end{definition}
Again, the relation \(\asequal{K_U}\) is an equivalence relation between both events and \(\sigma\)-algebras. 

Returning to counterfactual causal spaces, we can give precise mathematical definitions for the two extremes of how much information is shared between factual and counterfactual worlds, after an intervention. Let \(\fact{\sF}\subseteq\fact{\sH}\) be a factual \(\sigma\)-algebra and \(\cfact{\sF}\subseteq\cfact{\sH}\) a counterfactual \(\sigma\)-algebra. 
\begin{enumerate}
    \item If \(\fact{\sF}\indep{K_U}\cfact{\sF}\), then there is no shared information between \(\fact{\sF}\) and \(\cfact{\sF}\) after intervention on \(\sH_U\). 
    \item If \(\fact{\sF}\asequal{K_U}\cfact{\sF}\), then the information share between \(\fact{\sF}\) and \(\cfact{\sF}\) after intervention on \(\sH_U\) is maximal. 
\end{enumerate}
Instead of looking at sub-\(\sigma\)-algebras \(\fact{\sF}\) and \(\cfact{\sF}\) in the two worlds, we can also say that there is no shared information between the entire worlds after intervention on \(\sH_U\) if \(\fact{\sH}\indep{K_U}\cfact{\sH}\). However, if \(U\cap\fact{T}\neq\emptyset\) and \(U\cap\cfact{T}\neq\emptyset\), then it is not possible to have \(\fact{\sH}\asequal{K_U}\cfact{\sH}\), since, for any \(A\in\sH_{U\cap\fact{T}}\) and \(B\in\sH_{U\cap\cfact{T}}\), the interventional determinism axiom (\Cref{def:counterfactual_causal_space}\ref{item:interventional_determinism}) gives
\begin{alignat*}{2}
    K_U(\omega,A\cap B)=\ind_{A\cap B}(\omega)&\neq\ind_A(\omega)=K_U(\omega,A)\\
    &\neq\ind_B(\omega)=K_U(\omega,B),
\end{alignat*}
unless \(A=\Omega\) or \(B=\Omega\). This is the opposite of causal independence, since, for any \(A\in\sH_{U\cap \fact{T}}\) and \(B\in\sH_{U\cap \cfact{T}}\), we have, by the interventional determinism axiom again,
\begin{equation*}
    K_U(\omega,A\cap B)=\ind_{A\cap B}(\omega)=\ind_A(\omega)\ind_B(\omega)=K_U(\omega,A)K_U(\omega,B),
\end{equation*}
so \(A\) and \(B\) are always causally independent. 

The following result is about how causal independence translates to (conditional) independence after an intervention. The proofs are provided in \Cref{sec:proofs}. 
\begin{saveprop}{prop:independence_preserved}
    Let \(\cC=(\Omega,\sH,\PP,\KK)\) be a counterfactual causal space (\Cref{def:counterfactual_causal_space}), and let \(\cC^{\textnormal{do}(U,\QQ)}=(\Omega,\sH,\iPP{U}{\QQ},\iKK{U}{\QQ})\) be the counterfactual causal space obtained after intervening on \(\sH_U\) with a measure \(\QQ\) (\Cref{def:intervention_counterfactual}). 
    \begin{enumerate}[(i)]
        \item Let \(A,B\in\sH\) be events. If \(A\indep{K_U}B\) in \(\cC\) (see \Cref{def:causal_independence}), then \(A\) and \(B\) are conditionally independent given \(\sH_U\) in \(\cC^{\textnormal{do}(U,\QQ)}\). 
        \item If \(\sH_{U\cap \fact{T}}\indep{\QQ}\sH_{U\cap \cfact{T}}\) and \(\fact{\sH}\indep{K_U}\cfact{\sH}\) in \(\cC\), then \(\fact{\sH}\indep{\iPP{U}{\QQ}}\cfact{\sH}\). 
    \end{enumerate}
\end{saveprop}
In words, if \(A\) and \(B\) are causally independent on \(\sH_U\) in the original causal space, then they are conditionally independent given \(\sH_U\) after an intervention on \(\sH_U\). As an immediate corollary of \Cref{prop:independence_preserved}(i), for two sub-\(\sigma\)-algebras \(\sF_1,\sF_2\subseteq\sH\), if \(\sF_1\indep{K_U}\sF_2\), then for any measure \(\QQ\) on \(\sH_U\), they are conditionally independent under the measure \(\iPP{U}{\QQ}\) given \(\sH_U\). Further, \Cref{prop:independence_preserved}(ii) tells us that if the factual and counterfactual worlds are causally independent, and we intervene with a measure \(\QQ\) on \((\Omega,\sH_U)\) under which the worlds are independent, then in the resulting space, the worlds are (unconditionally) independent. 

The next result is about how causal synchronisation translates to synchronisation after an intervention. The proof is again provided in \Cref{sec:proofs}. 
\begin{saveprop}{prop:synchronisation_preserved}
    Let \(\cC=(\Omega,\sH,\PP,\KK)\) be a counterfactual causal space (\Cref{def:counterfactual_causal_space}), and let \(\cC^{\textnormal{do}(U,\QQ)}=(\Omega,\sH,\iPP{U}{\QQ},\iKK{U}{\QQ})\) be the counterfactual causal space obtained after intervening on \(\sH_U\) with a measure \(\QQ\) on \((\Omega,\sH_U)\) (\Cref{def:intervention_counterfactual}). Let \(A,B\in\sH\) be events. If \(A\asequal{K_U}B\) in \(\cC\), then \(A\asequal{\iPP{U}{\QQ}}B\) in \(\cC^{\textnormal{do}(U,\QQ)}\). 
\end{saveprop}
In words, if \(A\) and \(B\) are causally equal on \(\sH_U\) in the original causal space, then they are almost surely equal after the corresponding intervention. As an immediate corollary, for two sub-\(\sigma\)-algebras \(\sF_1,\sF_2\subseteq\sH\), if they are causally synchronised on \(\sH_U\) (\(\sF_1\asequal{K_U}\sF_2\)), then after an intervention on \(\sH_U\), they are synchronised (\(\sF_1\asequal{\iPP{U}{\QQ}}\sF_2\)).

\section{Multiple counterfactual worlds}\label{sec:multiple}
In \Cref{sec:counterfactual_probability_spaces,sec:counterfactual_causal_spaces}, we defined counterfactual spaces for two worlds. The aim of this section is to generalise to multiple worlds. 

\subsection{N-way counterfactual probability spaces}\label{subsec:N-way_counterfactual_probability_spaces}
We first define \(N\)-way counterfactual probability spaces. 

We take \(N\) sets of outcomes \(\Omega^1,...,\Omega^N\), and we equip each \(\Omega^j\), \(j=1,...,N\), with a \(\sigma\)-algebra \(\sE^j\). Then the entire measurable space \((\Omega,\sH)\) is obtained by taking the product of all the measurable spaces as follows:
\begin{equation*}
    (\Omega,\sH)=(\times_{j=1}^N\Omega^j,\otimes^N_{j=1}\sE^j).
\end{equation*}
For any outcome \(\omega\in\Omega\), we denote by \(\omega^j\) the projection of \(\omega\) to \(\Omega^j\), so that \(\omega\) is decomposed as \(\omega=(\omega^1,...,\omega^N)\). Also, for each \(j=1,...,N\), we denote by \(\sH^j\) the sub-\(\sigma\)-algebra of \(\sH\) consisting of cylinder sets \(\Omega^1\times...\times\Omega^{j-1}\times A\times\Omega^{j+1}\times...\times\Omega^N\). 

We are ready to define \(N\)-way counterfactual probability spaces. 
\begin{definition}\label{def:Nway_counterfactual_probability_space}
    An \emph{\(N\)-way counterfactual probability space} is defined as the triple \((\Omega,\sH,\PP)\), where \((\Omega,\sH)\) is a measurable space with the above product structure and \(\PP\) is a probability measure on \((\Omega,\sH)\). 
\end{definition}
Again, mathematically speaking, \(N\)-way counterfactual probability spaces are simply probability spaces with the above product structure on the measurable space \((\Omega,\sH)\). 1-way counterfactual probability spaces are simply probability spaces with no restrictions, and 2-way counterfactual probability spaces are precisely what was defined in \Cref{def:counterfactual_probability_space} in \Cref{sec:counterfactual_probability_spaces}. The specification of the measure \(\PP\) on events that do not live in a single \(\sH^j\) tells us how much information is shared between the worlds. 

\subsection{N-way counterfactual causal spaces}\label{subsec:n-way_counterfactual_causal_spaces}
We now generalise counterfactual causal spaces to \(N\) worlds. 

We take \(N\) index sets \(T^1,\dots,T^N\), and let \(T=\cup_{j=1,...,N}T^j\). For each \(j=1,...,N\) and each \(t\in T^j\), we take a set of outcome \(\Omega^j_t\) and equip it with a \(\sigma\)-algebra \(\sE^j_t\) (recall from footnote~\ref{foot:scripts} that we use superscripts for worlds and subscripts for components of worlds). For each \(j=1,...,N\), we define \(\Omega^j=\times_{t\in T^j}\Omega^j_t\) as the outcome set in the \(j^\text{th}\) world, and \(\sE^j=\otimes_{t\in T^j}\sE^j_t\) the corresponding \(\sigma\)-algebra. Then the entire measurable space \((\Omega,\sH)\) is obtained by taking the product of all the measurable spaces as follows:
\begin{equation*}
    (\Omega,\sH)=(\times_{j=1,...,N}\Omega^j,\otimes_{j=1,...,N}\sE^j),
\end{equation*}
and for any outcome \(\omega\in\Omega\), we denote by \(\omega^j\) the projection of \(\omega\) to \(\Omega^j\), so that the outcome \(\omega\) is decomposed as \(\omega=(\omega^1,...,\omega^N)\). Each \(\omega^j\) can be further decomposed as \(\omega^j=(\omega^j_t)_{t\in T^j}\), where, for each \(t\in T^j\), we denoted by \(\omega^j_t\) the projection of \(\omega^j\) onto \(\Omega^j_t\). For any \(S\in\sP(T)\), we denote by \(\sH_S\) the sub-\(\sigma\)-algebra of \(\sH\) generated by measurable rectangles \(\times_{j=1,...,N}(\times_{t\in T^j}A^j_t)\), where \(A^j_t\in\sE^j_t\) differ from \(\Omega^j_t\) only for finitely many \(t\) such that \(t\in S\). As a shorthand, for each \(j=1,...,N\), we write \(\sH^j=\sH_{T^j}\). 

We are ready to define an \(N\)-way counterfactual causal spaces. 
\begin{definition}\label{def:Nway_counterfactual_causal_space}
    An \emph{\(N\)-way counterfactual causal space} is defined as the quadruple \((\Omega,\sH,\PP,\KK)\), where \((\Omega,\sH)\) is a measurable space with the above product structure, \(\PP\) is a probability measure on \((\Omega,\sH)\) and \(\KK=\{K_S:S\subseteq T\}\), called the \emph{causal mechanism}, is a collection of transition probability kernels \(K_S\) from \((\Omega,\sH_S)\) into \((\Omega,\sH)\), called the \emph{causal kernel on \(\sH_S\)}, satisfying the following three axioms:
    \begin{enumerate}[(i)]
        \item\label{item:Nway_trivial_intervention} for all \(\omega\in\Omega\) and \(A\in\sH\), we have
        \begin{equation*}
            K_\emptyset(\omega,A)=\PP(A);
        \end{equation*}
        \item\label{item:Nway_no_cross_world_effect} for each \(j=1,...,N\), all \(\omega\in\Omega\), all \(\omega\in\Omega\), all \(S\in\sP(T)\) and all \(A\in\sH^j\), we have
        \begin{equation*}
            K_S(\omega,A)=K_{S\cap T^j}(\omega,A);
        \end{equation*}
        \item\label{item:Nway_interventional_determinism} for all \(A\in\sH_S\) and \(B\in\sH\), we have
        \begin{equation*}
            K_S(\omega,A\cap B)=\ind_A(\omega)K_S(\omega,B)=\delta_\omega(A)K_S(\omega,B);
        \end{equation*}
        in particular, for \(A\in\sH_S\), we have
        \begin{equation*}
            K_S(\omega,A)=\ind_A(\omega).
        \end{equation*}
    \end{enumerate}
\end{definition}
It should be remarked again that \(N\)-way counterfactual causal spaces are special cases of causal spaces with the additional axiom of no cross-world causal effect. Clearly, \(1\)-way counterfactual causal spaces are simply causal spaces with no other restrictions than the usual axioms of causal spaces, and \(2\)-way counterfactual causal spaces are precisely what was defined in \Cref{def:counterfactual_causal_space} in \Cref{sec:counterfactual_causal_spaces}. 

For the sake of completeness, we define interventions in \(N\)-way counterfactual causal spaces, in an analogous way to \Cref{def:intervention_counterfactual}. 
\begin{definition}\label{def:intervention_Nway_counterfactual}
    Let us take an \(N\)-way counterfactual causal space \((\Omega,\sH,\PP,\KK)\) as in \Cref{def:Nway_counterfactual_causal_space}, a subset \(U\subseteq T\) and a probability measure \(\QQ\) on \((\Omega,\sH_U)\). An \emph{intervention on \(\sH_U\) via \(\QQ\)} yields a new \(N\)-way counterfactual causal space
    \begin{equation*}
        (\Omega,\sH,\iPP{U}{\QQ},\iKK{U}{\QQ}),
    \end{equation*}
    where the \emph{intervention measure} \(\iPP{U}{\QQ}\) is a probability measure on \((\Omega,\sH)\) defined, for \(A\in\sH\), by
    \begin{equation*}
        \iPP{U}{\QQ}(A)=\int\QQ(d\omega)K_U(\omega,A)
    \end{equation*}
    and \(\iKK{U}{\QQ}=\{\iK{U}{\QQ}{S}:S\subseteq T\}\) is the \emph{intervention causal mechanism} whose \emph{intervention causal kernels} are
    \begin{equation*}
        \iK{U}{\QQ}{S}(\omega_S,A)=\int\QQ(d\omega'_{U\setminus S})K_{S\cup U}((\omega_{S},\omega'_{U\setminus S}),A).
    \end{equation*}
\end{definition}
It must be checked that the \(N\)-way counterfactual causal space obtained after an intervention is indeed an \(N\)-way counterfactual causal space, i.e. the intervention causal mechanism \(\iKK{U}{\QQ}\) satisfies the axioms of \Cref{def:Nway_counterfactual_causal_space}. The following theorem proves this, and in doing so, proves \Cref{thm:intervention} as a special case of \(2\)-way counterfactual causal spaces. The proof is given in \Cref{sec:proofs}. 
\begin{savetheorem}{thm:Nway_interventions}
    The intervention causal mechanism \(\iKK{U}{\QQ}\) given in \Cref{def:intervention_Nway_counterfactual} satisfies the axioms of causal mechanisms given in \Cref{def:Nway_counterfactual_causal_space}. 
\end{savetheorem}

\section{Related Works}\label{sec:related_works}
Counterfactuals have been extensively studied by philosophers in the tradition of possible worlds semantics, with influential accounts given by \citet{goodman1947problem,stalnaker1968theory} and \citet{lewis1973counterfactuals}. Further, psychologists have studied the significant role that counterfactual thinking plays in a child's development, perception and reasoning in adults, and its impact on decision-making, emotions and biases \citep{byrne2016counterfactual,waldmann2017oxford}. Counterfactuals also feature prominently in a wide range of application areas, such as fairness \citep{kusner2017counterfactual,garg2019counterfactual,rosenblatt2023counterfactual}, harm \citep{richens2022counterfactual,beckers2022causal,beckers2023quantifying,straitouri2024controlling}, interpretable machine learning through counterfactual explanations and algorithmic recourse \citep{guidotti2024counterfactual,dissanayake2024model,verma2024counterfactual}, counterfactual image editing \citep{pan2024counterfactual}, and counterfactual sampling and generations \citep{ribeiro2023high,hao2024natural,melistas2024benchmarking,jung2024counterfactually,raghavan2024counterfactual}, with clinical applications \citep{degrave2023dissection,lee2024clinical}. Establishing a rigorous, axiomatic mathematical framework for counterfactuals is a crucial endeavour, laying the foundation for any kind of quantitative research involving estimation of, or reasoning with, counterfactual probabilities. 

In the rest of this section, our review of the related formalisms largely focus on the two major frameworks of counterfactuals (in fact, of causality) mentioned in the introduction (\Cref{sec:introduction}), namely, the SCMs and POs. We show that, starting from a specification of an SCM or a PO framework, we can construct a counterfactual space, demonstrating the fact that counterfactual spaces strictly generalise the existing formalisms. 

\begin{remark}
    In a series of papers, \citet{dawid1999needs,dawid2000causal,dawid2007counterfactuals} makes a clear distinction between \emph{effects of causes} and \emph{causes of effects}. He argues that counterfactual considerations are unnecessary and potentially misleading for effects of causes, and only interventional considerations are required (for which he proposes a decision-theoretic framework). He also argues that inferring causes of effects---which in turn requires thinking about counterfactuals---is impossible to corroborate with data, and that, since it is not suitably empirical, there is no point in developing a theory of it. We acknowledge that, without assumptions, real-world validation of counterfactuals is impossible, and also that counterfactuals are often irrelevant for causality, as our orthogonal view (\Cref{fig:ladder}) shows. However, we do not agree that this provides grounds for an outright rejection of a formalism for counterfactuals; even though empirical verification may be impossible, axiomatising such an fundamental component of human thought is still, for reasons we discuss, a worthwhile pursuit.
\end{remark}

\subsection{Structural causal models (SCMs)}\label{subsec:scms}
Pearl's SCMs \citep{pearl2009causality,peters2017elements} remain one of the most influential and widely used mathematical frameworks for counterfactuals, and for causality as a whole, with several variants to accommodate different desiderata \citep{hiddleston2005causal,rips2010two,fisher2017counterlegal,lee2017hiddleston,bongers2021foundations}. However, despite all its merits and appealing properties, this framework has some major, well-known limitations as a foundational axiomatisation, as discussed in the introduction (\Cref{sec:introduction}).  
Thus, though we submit 
that SCMs provide a valuable tool to treat a specific type of counterfactuals, we object to the assertion of \citet{pearl2000causal} that 
\begin{quote}
    \say{Functional models, in the form of nonparametric structural equations, thus provide both formal semantics and conceptual basis for a complete mathematical theory of counterfactuals}.
\end{quote}
Let us now recall the mathematical definition of an SCM. An SCM is a triple \(\cM=(\bU,\bV,\bF)\), where \(\bU=\{U_1,...,U_m\}\) is a set of exogenous variables, \(\bV=\{V_1,...,V_n\}\) is a set of endogenous variables with each \(V_i\) taking values in the measurable space \((\Omega_i,\sE_i)\) for \(i=1,...,n\), and \(\bF=\{f_1,...,f_n\}\) are the structural equations such that \(V_i=f_i(\bPA_i,\bU_i)\) for \(i=1,...,n\), with \(\bPA_i\subseteq\bV\setminus\{V_i\}\) and \(\bU_i\subseteq\bU\). Hence, any subset of the endogenous variables \(\bX\subseteq\bV\) is a deterministic function of the exogenous variables \(\bU\). Given a specific value \(\bu\) of \(\bU\), we write \(\bX(\bu)\) for the value of \(\bX\) determined by \(\bu\). 

We make the model probabilistic by imposing a measure \(\PP^\bU\) on \(\bU\), which induces a measure on \(\bV\) as a pushforward measure. Specifically, for an event \(A\in\otimes_{i=1}^n\sE_i\),
\begin{equation*}
    \PP(A)=\int\PP^\bU(d\bu)\ind\{\bV(\bu)\in A\}.
\end{equation*}
With a slight abuse of notation, for a subset \(\bX\) of \(\bV\), we write \(\Omega_{\bX}=\times_{i\in[n],V_i\in\bX}\Omega_i\) and \(\sE_{\bX}=\otimes_{i\in[n],V_i\in\bX}\sE_i\). For a realisation \(\bx\) of \(\bX\), the \emph{sub-model} \(\cM_{\bX=\bx}=(\bU,\bV,\bF_{\bX=\bx})\) of \(\cM\) is given by \(\bF_{\bX=\bx}=\{f_i:V_i\notin\bX\}\cup\{\bX=\bx\}\), and the \emph{potential response} of \(\bY\subseteq\bV\) to the action \(\text{do}(\bX=\bx)\) under the noise values \(\bu\) is denoted as \(\bY_{\bX=\bx}(\bu)\in\Omega_{\bY}\). 

Using this model, the \emph{probability of counterfactuals} are calculated as follows. Let us consider two identical SCMs \(\cM=(\bU,\bV,\bF)\) and \(\cM^*=(\bU^*,\bV^*,\bF^*)\). For any sets of variables \(\bY^*,\bX^*\subseteq\bV^*\) and \(\bZ,\bW\subseteq\bV\) and events \(A\in\sE_{\bY^*}\) and \(B\in\sE_\bZ\), we have
\begin{equation*}
    \PP(\bY^*_{\bX^*=\bx^*}\in A,\bZ_{\bW=\bw}\in B)=\int\PP^\bU(d\bu)\ind\{\bY^*_{\bX^*=\bx^*}(\bu)\in A\}\ind\{\bZ_{\bW=\bw}(\bu)\in B\},
\end{equation*}
where \(\bY^*_{\bX^*=\bx^*}\) and \(\bZ_{\bW=\bw}\) are potential responses from sub-models \(\cM_{\bX^*=\bx^*}^*\) and \(\cM_{\bW=\bw}\) that share the same values of the exogenous variables \(\bU\). 

In this framework, the type of counterfactuals most commonly considered is of the form \(\PP_{\sigma\bZ}(\bY^*_{\bX^*=\bx^*}\in A)\), where \(\sigma\bZ\) 
is the \(\sigma\)-algebra generated by the random variable \(\bZ\) (the so-called abduction--action--prediction procedure). The intervention \(\bX^*=\bx^*\) and the observation \(\bZ\) may be incompatible. This is performed simply by taking the conditional distribution \(\PP^\bU_{\sigma\bZ}(\cdot)\) given \(\sigma\bZ\) and the sub-model \(\cM_{\bX^*=\bx^*}\):
\begin{equation*}
    \PP_{\sigma\bZ}(\bY^*_{\bX^*=\bx^*}\in A)=\int\PP^\bU_{\sigma\bZ}(d\bu)\ind\{\bY^*_{\bX^*=\bx^*}(\bu)\in A\}.
\end{equation*}
Suppose that we have an arbitrary specification of an SCM as given above. Then we specify a counterfactual causal space \((\Omega,\sH,\PP,\KK)\) as follows.
\begin{itemize}
    \item We let \(\fact{T}=\cfact{T}=[n]\), and take \(T=\fact{T}\cup \cfact{T}\). We also let \(\fact{\Omega}=\cfact{\Omega}=\times_{i=1}^n\Omega_i\), and \(\Omega=\fact{\Omega}\times\cfact{\Omega}\). Finally, we let \(\fact{\sE}=\cfact{\sE}=\otimes_{i=1}^n\sE_i\), and \(\sH=\fact{\sE}\otimes\cfact{\sE}\). 
    \item For any measurable rectangle \(A\times B\in\sH\) with \(A\in\fact{\sE}\) and \(B\in\cfact{\sE}\), we have
    \begin{equation*}
        \PP(A\times B)=\int\PP^\bU(d\bu)\ind\{\bV(\bu)\in A\}\ind\{\bV^*(\bu)\in B\}.
    \end{equation*}
    This is extended to all of \(\sH\) in the usual way. 
    \item Take any \(S\in\sP(T)\), and write \(\bX=\{V_i\in\bV:i\in\fact{T}\cap S\}\) and \(\bX^*=\{V_i\in\bV^*:i\in\cfact{T}\cap S\}\) for the variables being intervened on in the factual and counterfactual worlds respectively. Then the corresponding causal kernel for a rectangle \(A\times B\in\sH\) with \(A\in\fact{\sE}\) and \(B\in\cfact{\sE}\) is given by
    \begin{equation*}
        K_S((\bx,\bx^*),A\times B)=\int\PP^\bU(d\bu)\ind\{\bV_{\bx}(\bu)\in A\}\ind\{\bV^*_{\bX^*=\bx^*}(\bu)\in B\}.
    \end{equation*}
    In other words, we are taking the pushforward measure of \(\PP^\bU\) through the new structural equations in the sub-models \(\cM_{\bx}\) and \(\cM^*_{\bx^*}\) given by the interventions. It is again extended to all events in \(\sH\) in the usual way. One should think of this as the causal kernel corresponding to the interventions \(\bX=\bx\) and \(\bX^*=\bx^*\).
\end{itemize}
We see that an SCM uniquely determines a corresponding counterfactual causal space, and so counterfactual causal spaces generalise all of observational, interventional and counterfactual information of SCMs. 

The standard SCMs discussed above force the pre-intervention worlds to be synchronised, since the exogenous variables and the structural equations are shared. To relax this constraint, \citet{von2023backtracking} introduced \emph{backtracking} SCMs, which allow a more general level of information share between the worlds. We review this formalism below, and for each specification of a backtracking SCM, construct counterfactual probability space that has the same counterfactual information. 

Let us take two identical SCMs \(\cM=(\bU,\bV,\bF)\) and \(\cM^*=(\bU^*,\bV^*,\bF^*)\), for the factual and counterfactual worlds respectively. Backtracking SCMs define a backtracking measure \(\PP^\text{B}\) over the exogenous noise variables \(\bU\) and \(\bU^*\). Then for events \(A,B\in\otimes^n_{i=1}\sE_i\), the backtracking counterfactual probabilities are given by
\begin{equation*}
    \PP^\textnormal{B}(\bV\in A,\bV^*\in B)=\int\ind\{\bV(\bu)\in A\}\ind\{\bV^*(\bu^*)\in B\}\PP^\text{B}(d\bu,d\bu^*).
\end{equation*}
No intervention takes place in either \(\cM\) or \(\cM^*\). 

To construct the corresponding counterfactual probability space, we first construct a measurable space \((\Omega,\sH)\) by letting \(\fact{\Omega}=\cfact{\Omega}=\times_{i=1}^n\Omega_i\) and \(\Omega=\fact{\Omega}\times\cfact{\Omega}\), and \(\fact{\sE}=\cfact{\sE}=\otimes_{i=1}^n\mathscr{E}_i\) and \(\sH=\fact{\sE}\otimes\cfact{\sE}\). Then, we define \(\PP\) on the rectangles \(A\times B\) for \(A\in\fact{\sE}\) and \(B\in\cfact{\sE}\) by 
\[\PP(A\times B)=\int\ind\{\bV(\bu)\in A\}\ind\{\bV^*(\bu^*)\in B\}\PP^\text{B}(d\bu,d\bu^*).\]
We extend \(\PP\) to all measurable sets in \(\sH\) in the usual way.

\subsection{Potential outcomes}\label{subsec:po}
A major competing framework of causality and counterfactuals is the \emph{potential outcomes} framework \citep{imbens2015causal,hernan2020what}, widely adopted in, for example, social and biomedical sciences, and econometrics. We argue that, while this framework has many well-established virtues, similar to SCMs, it also falls short as a foundational axiomatisation of counterfactuals. Some of its most obvious limitations are similar to those of SCMs, in that it struggles with cycles or continuous-time stochastic processes. Further, it is a \emph{static} framework in which no changes to the mathematical quantities (most notably, the probability distributions) can take place, despite the fact that the effect of interventions are arguably precisely such changes.
The framework simply adds \say{potential outcome variables} to the model, which represent what \emph{would} happen \emph{if} an intervention were to take place. As a result, no consideration of sequential interventions, for example, is built in. As for counterfactuals, only those that are based on different values of the treatment variable are incorporated; no consideration of non-interventional counterfactuals, or of stochastic interventions in at least one of the worlds, is possible. 

In the potential outcomes framework, most often, a treatment variable, an outcome variable and covariate variables are designated a priori. However, we adopt a more general definition of \citet[Definition 1]{ibeling2023comparing} (which, in turn, is based on \citet{holland1986statistics}, and is named the \say{Rubin causal model}). Here, we have a given set of endogenous variables, and any of these variables can act as the treatment or the outcome.

Suppose, as in the SCM framework, that \(\bV=\{V_1,...,V_n\}\) is a set of endogenous variables, with each \(V_i\) taking values in a measurable space \((\Omega_i,\mathcal{E}_i)\). We also have a finite set \(\cU\) of \emph{units}, and a probability measure \(\PP^\cU\). We take a set \(\cD\) of \say{potential outcome variables}, of the form \(V_{i,\bx}\) for some \(V_i\in\bV\), \(\bX\subseteq\bV\) and some value \(\bx\in\Omega_\bX\). The potential outcome \(V_{i,\bx}\) takes values in \(\Omega_i\). We finally have a set of functions \(\bF\) which consists of a function \(f_{i,\bx}:\cU\to\Omega_i\) for each \(V_{i,\bx}\in\cD\) and a function \(f_i:\cU\to\Omega_i\) for each \(i\in\{1,...,n\}\). The whole model is the quintuple \((\cU,\bV,\cD,\bF,\PP^\cU)\), and we get a joint measure over the endogenous variables and the potential outcomes by a pushforward of \(\PP^\cU\) through the functions in \(\bF\). 

Of course, if we only considered potential outcomes \(V_{i,x}\) for a single \(i\in\{1,\dots,n\}\) and a single variable \(X\in\bV\), then we would recover the common case in which the outcome and treatment variables are fixed in advance: \(V_i\) and \(X\) respectively.

We now show that we can uniquely construct an \(N\)-way counterfactual probability space (\Cref{def:Nway_counterfactual_probability_space}) from an arbitrary specification of the potential outcomes framework. The number of counterfactuals that are considered in the potential outcomes framework is not the number of counterfactual worlds of interest, but the number of treatment variable values of interest. In the most common case of binary treatment, we need a \(3\)-way counterfactual probability space, one world for the \say{observed} variables, and one each for the values of the treatment variable. 

Take a specification of the potential outcomes framework \((\cU,\bV,\cD,\bF,\PP^\cU)\). Then we take an \((N+1)\)-way counterfactual probability space, where \(N\) is the number of distinct values \(\bx\) in the subscript of the potential outcomes in \(\cD\), which we enumerate as \(\{\bx^1,...,\bx^N\}\). For \(j=N+1\), corresponding to the observed world, we take the measurable space \((\Omega^{N+1},\sE^{N+1})=\otimes_{i=1}^n(\Omega_i,\sE_i)\), the domain of the entire set of endogenous variables. For each \(j=1,...,N\), we take the measurable space \((\Omega^j,\sE^j)=\otimes_{i\in S^j}(\Omega_i,\sE_i)\), where \(S^j=\{i:V_{i,\bx^j}\in\cD\}\). Then the entire measurable space is given by \((\Omega,\sH)=\otimes^{N+1}_{j=1}(\Omega^j,\sE^j)\). 

The measure \(\PP\) on \((\Omega,\sH)\) is given as follows. For a rectangular event in \(\sH\) of the form \(\otimes_{j=1}^{N+1}\otimes_{i\in S^j}A^j_i\) with \(A^j_i\in\sE_i\) for each \(j\in\{1,...,N+1\}\), we define
\begin{equation*}
    \PP\left(\times_{j=1}^{N+1}\times_{i\in S^j}A^j_i\right)=\int\prod_{j=1}^{N+1}\prod_{i\in S^j}\ind\{f_{i,\bx^j}(\bu)\in A^j_i\}\PP^\cU(d\bu).
\end{equation*}
These rectangles generate \(\sH\), so we can extend \(\PP\) in the usual way to all of \(\sH\). 

In the potential outcomes framework, no new mathematical object is introduced to encode causality: it is simply read off from the single probability measure over all the variables, including the potential outcomes which represent what \emph{would} happen \emph{if} an intervention were to take place. No changes to the mathematical quantities, in particular on the measure, takes place. It is by reason that counterfactual \emph{probability} spaces, not counterfactual \emph{causal} spaces, were used in this section. By assigning the potential outcomes in the appropriate counterfactual worlds, we constructed the counterfactual probability spaces that corresponds exactly to given specifications of the potential outcomes framework. This stands in contrast to the SCM, causal space or counterfactual causal space frameworks, in which an intervention leads to a change of the measure. 

\section{Conclusion}\label{sec:conclusion}
In this paper, we introduced \emph{counterfactual probability spaces} and \emph{counterfactual causal spaces} as axiomatic frameworks for capturing counterfactuals, rigorously grounded in measure theory. They are special cases of probability spaces and causal spaces, which are respectively measure-theoretic axiomatisations of the concepts of probability and of interventions. We viewed interventional causality and counterfactuals as orthogonal concepts, which we brought together in a single framework of counterfactual causal spaces. 

We suggested that the essence behind the study of counterfactuals is the simultaneous consideration of two (or more) parallel worlds, and we proposed a way of capturing the shared information between the worlds. In counterfactual probability spaces, and in counterfactual causal spaces before intervention, the shared information is encoded in the probability measure, and after an intervention in counterfactual causal spaces, it is encoded in the corresponding causal kernel. The two extremes of the extent to which the worlds are related are captured in the definitions of \emph{independence} and \emph{synchronisation}; these possibilities are either impossible or imposed by definition in prominent frameworks. We have shown that our spaces strictly subsume all previous formalisms, while dispensing with major assumptions that are required in their definitions, such as acyclicity, discreteness, and that endogenous variables do not causally affect the exogenous variables. 

We demonstrated that the definitions of conditional causal effects have a natural interpretation in counterfactual causal spaces, and can be used to answer queries that commonly arise in ordinary human thought. 

As outlined in \Cref{sec:related_works}, counterfactuals have found a wealth of applications, and will no doubt continue to do so, especially with the advance of artificial intelligence and generative models \citep{geiger2025causal}. We believe that a rigorous, axiomatic treatment of this fundamental concept will lay a valuable foundation for future research endeavours involving counterfactuals. 

As a final note, it is not our intention to criticise existing frameworks of counterfactuals, nor to replace them. We believe that those already in the literature, such as the SCMs, potential outcomes or single-world intervention graphs (SWIGS; \citep{richardson2013single}, which combine the potential outcomes approach with graphical approaches\footnote{We did not explicitly review SWIGs in this paper, because we believe that their goal is not to provide a foundational framework of counterfactuals, but rather to provide a useful graphical tool for dealing with counterfactuals, in particular, Markov conditions representing conditional independencies among factual and counterfactual variables. The SCM and the potential outcomes framework are often claimed to play the role of a foundational framework. We hope to have cast doubt on such claims throughout this paper.}) are useful frameworks, that will no doubt continue to play an important role in the theory of counterfactuals. However, as we argued throughout this paper, they rely on assumptions \emph{by definition}, and/or are unable to represent certain kinds of counterfactual scenarios, and hence fall short as a foundational, axiomatic framework for counterfactuals. We believe that the existing frameworks will continue to play important roles in elegantly and succinctly \emph{specifying} a counterfactual space, just as (for example) various parametric distributions in probability theory play the role of specifying a probability space. 

\vspace{1cm}
\bibliographystyle{unsrtnat}
\bibliography{ref} 

\clearpage
\begin{appendix}
\section{Actual causality}\label{sec:actual_causality}
There is an active area of research called \emph{actual causality} (or \emph{token causality}) \citep{halpern2001causes,halpern2005causes,halpern2015modification,halpern2016actual,beckers2018principled,beckers2021causal}. Here, one asks, after having observed A and B, queries of the form \say{Did A cause B in the specific situation that played out in the factual world?} (e.g. did my father's smoking habit of 30 years cause his lung cancer?). Actual causality differs from the notion of \emph{general causality} (or \emph{type causality}), in which one is interested in causal questions about a general phenomenon (does smoking increase one's chances of getting lung cancer in general?). The investigation of actual causality is especially relevant in the realm of law, where blame must be assigned to individuals for particular occurrences of events \citep{halpern2015cause,fischer2024broken}, but also to highlight suitable targets for intervention \citep{hitchcock2009cause}. The general theory of causal spaces, including the theory of counterfactual causal spaces presented in the main body of this paper, is concerned with general, rather than actual, causality. 

In contrast to general causality, which has been precisely defined in all existing frameworks of causality, a universal mathematical definition of actual causality has been elusive \citep{beckers2021causal,fischer2024actual}. In the past couple of decades, many proposals were made for this purpose (predominantly in the SCM framework), counterexamples were found, new proposals were made, and the process is on-going \citep{halpern2001causes,halpern2005causes,halpern2015modification,halpern2015graded,halpern2016actual,icard2017normality,beckers2021causal,gallow2021model,andreas2021ramsey}. In fact, some authors argue for a pluralist approach to actual causality, instead of searching for \emph{the} definition of actual causality \citep{fischer2024actual,fischer2024three}. 

Once we narrow down the definition of causal effects to \emph{conditional effects of individual outcomes}, the queries are of the form, \say{does a specific outcome have a causal effect in the specific situation observed in the factual world?} This, at first glance, sounds similar to actual causality. However, a closer investigation reveals that there are subtle but irreconcilable differences, both philosophical and mathematical, between what we can answer with the definition of conditional causal effects in causal spaces and questions asked in actual causality. 

In this section, we first review a subset of the definitions of actual causality proposed in the SCM framework. These can, of course, be translated into counterfactual causal spaces, since we have already shown that counterfactual causal spaces strictly generalise SCMs. However, we will argue, with the running example from the main body, why these existing definitions should not, in our opinion, be \emph{the} definition of actual causality. 

In a series of papers \citep{halpern2001causes,halpern2005causes,halpern2015modification}, Halpern and Pearl proposed mathematical definitions of actual causality, and they remain the most influential account of actual causality (see also the book, \citep{halpern2016actual}). In a series of papers \citep{beckers2018principled,beckers2021counterfactual,beckers2021causal,beckers2025actual}, Beckers also studied the problem of actual causality within the SCM framework. Both lines of work focus largely on the discrete case (where an \say{event} is synonymous with a variable taking a particular value, rather than the definition of events in probability theory) and the discussion of probabilities is largely side-stepped. The underlying philosophy is that \say{A is an actual cause of B if A is a \emph{necessary element of a sufficient set}} (NESS). Here, we review a subset of those definitions. 

Recall that an SCM is a triple, \(\cM=(\bU,\bV,\bF)\), where \(\bU\) is a set of exogenous variables. In the theory of Halpern and Pearl, actual causality is defined for a particular outcome value \(\bu\) of the exogenous variables \(\bU\). Recall also that, for a set of endogenous variables \(\bX\subseteq\bV\), we write \(\bX(\bu)\) for the particular realisation of \(\bX\) (deterministically) induced by the noise value \(\bu\) through the structural equations \(\bF\). 

\begin{definition}\label{def:actual_causality}
    We say that \(\bX=\bx\) is an \emph{actual cause} of \(Y=y\), if the following three conditions hold. 
    \begin{description}
        \item[AC1] \(\bX(\bu)=\bx\) and \(Y(\bu)=y\).
        \item[AC2] See below.
        \item[AC3] \(\bX\) is minimal; there is no strict subset \(\bX'\) of \(\bX\) such that \(\bX'=\bx'\) satisfies conditions AC1 and AC2, where \(\bx'\) is the restriction of \(\bx\) to the variables in \(\bX'\).
    \end{description}
\end{definition}
Many proposals of actual cause over the years have kept this format, where the conditions AC1 and AC3 remain the same. We review four of the most prominent proposals for AC2, given respectively in \citep{halpern2001causes,halpern2005causes,halpern2015modification,beckers2021causal}. Recall from \Cref{subsec:scms} that the potential response of a variable \(Y\) after the intervention \(\bX=\bx\) with the noise value \(\bu\) is written as \(Y_{\bX=\bx}(\bu)\). 
\begin{itemize}
    \item \textbf{Original HP definition}
    \begin{description}
        \item[AC2(necessity)] There is a partition of \(\bV\) into two disjoint subsets \(\bZ\) and \(\bW\) with \(\bX\subseteq\bZ\) and a setting \(\bx'\) and \(\bw\) of the variables in \(\bX\) and \(\bW\), respectively, such that \(Y_{\bX=\bx',\bW=\bw}\neq y\).
        \item[AC2(sufficiency)] If the value \(\bz^*\) is such that \(\bZ(\bu)=\bz^*\), then for all subsets \(\bZ'\) of \(\bZ\setminus\bX\), we have \(Y_{\bX=\bx,\bW=\bw,\bZ^*=\bz^*}(\bu)=y\). 
    \end{description}
    \item \textbf{Updated HP definition}
    \begin{description}
        \item[AC2(necessity)] There is a partition of \(\bV\) into two disjoint subsets \(\bZ\) and \(\bW\) with \(\bX\subseteq\bZ\) and a setting \(\bx'\) and \(\bw\) of the variables in \(\bX\) and \(\bW\), respectively, such that \(Y_{\bX=\bx',\bW=\bw}\neq y\).
        \item[AC2(sufficiency)] If \(\bz^*\) is such that \(\bZ(\bu)=\bz^*\), then for all subsets \(\bW'\) of \(\bW\) and all subsets \(\bZ'\) of \(\bZ\setminus\bX\), we have \(Y_{\bX=\bx,\bW'=\bw,\bZ^*=\bz^*}(\bu)=y\). 
    \end{description}
    \item \textbf{Modified HP definition}
    \begin{description}
        \item[AC2] There is a set \(\bW\) of variables in \(\bV\) and a setting \(\bx'\) of the variables in \(\bX\) such that, if \(\bW(\bu)=\bw^*\), then \(Y_{\bX=\bx',\bW=\bw^*}(\bu)\neq y\). 
    \end{description}
    \item \textbf{Def 2 of Beckers}
    \begin{description}
        \item[AC2(necessity)] There exist sets \(\bW\) and \(\bN\) with \(Y\in\bN\), and values \(\bx'\), such that for all \(\bS\subseteq\bN\) with \(Y\in\bS\), and for all \(\bs\in\Omega_\bS\) such that \(y\in\bs\), there exists a \(\bt\in\Omega_{\bV\setminus\{\bX\cup\bW\cup\bS\}}\) such that \(\bS_{\bX=\bx',\bW=\bw^*,\bT=\bt}(\bu)\neq\bs\), where \(\bw^*\) is such that \(\bW(\bu)=\bw^*\).
        \item[AC2(sufficiency)] For all \(\bc\in\Omega_{\bV\setminus\{\bX\cup\bW\cup\bN\}}\), we have that \(\bN_{\bX=\bx,\bW=\bw^*,\bC=\bc}(\bu)=\bn^*\), where \(\bw^*,\bn^*\) are such that \(\bW(\bu)=\bw^*\) and \(\bN(\bu)=\bn^*\). 
    \end{description}
\end{itemize}
Without going into the details, we note that all of these definitions rely two crucial assumptions:
\begin{enumerate}
    \item that the worlds are synchronised under the same value \(\bu\) of the exogenous variables \(\bU\), both in the observational state and after any intervention;
    \item that the observational and interventional distributions are coupled through the structural equations \(\bF\). 
\end{enumerate}
Hence, in all situations where this assumption is violated, these definitions lose the grounds they stand on. In general, actual causality is a problem of inferring \emph{causality} from \emph{observations}. Hence, in full generality, it is an ill-posed problem, even if one has access to a counterfactual world in which we can perform interventions, without imposing assumptions on how the observational and interventional distributions are related. The fact that a principled definition of actual causality has been so elusive suggests that the assumptions imposed by the SCM framework are not sufficient for the purpose. It is beyond the scope of this paper to propose a set of assumptions that would accommodate a definition. 

We give an example of a case which is not catered for by the SCM framework, and therefore the definitions of actual causality therein, by continuing our running example, \Cref{ex:exam}. 

\begin{example}\label{ex:exam_cycle}
    In addition to the observational distribution and causal kernels corresponding to intervening on \(\class\) given in \Cref{tab:exam_observation,tab:exam_class,tab:exam_noclass}, we define the causal kernel of intervening on the exam result in the counterfactual world in \Cref{tab:exam_pass}. 
    \begin{table}[t]
        \centering
        \begin{tabular}{cccccc|c}
            &&\multicolumn{5}{c}{Counterfactual}\\
            &\(K_{\cfact{\exam}}(P,\cdot)\)&\((Y,P)\)&\((Y,F)\)&\((N,P)\)&\((N,F)\) & Sum \\
            \multirow{5}{*}{Factual }&\((Y,P)\)&0&0&0.43&0&0.43\\
            &\((Y,F)\)&0&0&0.21&0&0.21\\
            &\((N,P)\)&0&0&0.19&0&0.19\\
            &\((N,F)\)&0&0&0.17&0&0.17\\\cline{2-7}
            &\text{Sum}&0&0&1&0&1\\
        \end{tabular}
        \caption{The causal kernel for intervening on the exam result to be a pass in the counterfactual world.}
        \label{tab:exam_pass}
    \end{table}

    The causal kernel reflects the fact that, if the students are told that they will receive a pass grade no matter what, then no student will attend the revision class. We make some remarks on this causal kernel. 
    \begin{enumerate}[(i)]
        \item The no cross-world causal effect and interventional determinism axioms (\Cref{def:counterfactual_causal_space}) hold. 
        \item We have a cyclic causal relationship---the attendance at the revision class causally affects (albeit weakly) the exam results, as seen in \Cref{tab:exam_class,tab:exam_noclass}, and the intervention on the exam result has a very strong causal effect on the attendance at the revision class. All causal kernels defined in \Cref{tab:exam_class,tab:exam_noclass,tab:exam_pass} are valid, and constitute parts of a valid counterfactual causal space. However, no SCM can be constructed with these observational and interventional distributions. 
        \item Since the observational distribution (\Cref{tab:exam_observation}) tells us nothing about what would happen under the intervention given in \Cref{tab:exam_pass}, being able to intervene in the counterfactual world tells us nothing about whether the exam result had an \say{actual causal effect} in the observational state. 
    \end{enumerate}
\end{example}
As mentioned earlier, there is a difference between actual causality and what the definition of conditional causal effect allows us to answer by conditioning on the factual world and intervening on the counterfactual world. 
\begin{description}
    \item[Conditional causal effect] Conditioning on the factual world refines the situation under investigation in the counterfactual world, depending on what information is transmitted across worlds. Given this refined situation, conditional causal effect asks, \say{what would (have) happen(ed) if we intervened in this situation?} 
    \item[Actual causality] The question asked in actual causality is fundamentally different: it is of the form, \say{in the observed situation, what was the cause?} 
\end{description}
Existing works try to answer actual causality by intervening in the counterfactual world. But the intervention carried out does not correspond to the query interested in, instead answering one of the form, \say{given that we observed \(\{\bx,y\}\), would we still have had \(y\) if we intervened with \(\bx'\) instead?} This is subtly but unquestionably different to the above question, asked in actual causality. While the existing definitions (\Cref{def:actual_causality}) are valuable efforts to bridge the two types of queries, the fundamental block is that we are trying to infer a causal effect from an observation. For this, we need assumptions on how the interventional distribution corresponding to the causal effect of interest is connected to the observational distribution, which are not afforded by the assumptions of the SCM framework. 

We leave it as future work to study assumptions that will facilitate the study of actual causality in counterfactual causal spaces. 

\section{Sources}\label{sec:sources}
We begin this section by recalling the definition of \emph{sources}, which, even though it was originally defined only for causal spaces, extends to counterfactual causal spaces without any adjustments. Sources are those \(\sigma\)-algebras on which the causal kernel and the conditional probability coincide almost surely (with respect to the observational measure \(\PP\)). 
\begin{definition}\label{def:sources}[{\citet[Definition D.1]{park2023measure}}]
    Let \((\Omega,\sH,\PP,\KK)\) be a causal space as defined in \Cref{def:causal_space}. Let \(U\subseteq T\) be a subset, \(A\in\sH\) an event and \(\sF\) a sub-\(\sigma\)-algebra of \(\sH\). We say that
    \begin{enumerate}
        \item \(\sH_U\) is a \emph{(local) source} of \(A\) if, for \(\PP\)-almost all \(\omega\in\Omega\), we have \(K_U(\omega,A)=\PP_{\sH_U}(\omega,A)\);
        \item \(\sH_U\) is a \emph{(local) source} of \(\sF\) if \(\sH_U\) is a source of all \(A\in\sF\); and
        \item \(\sH_U\) is a \emph{global source} if \(\sH_U\) is a source of all \(A\in\sH\).
    \end{enumerate}
\end{definition}
In the causal inference community, there is a very strong focus on the problem of \emph{identifiability}, i.e. the problem of inferring causal information using just the observational data. In terms of (counterfactual) causal spaces, it is the problem of inferring information about the causal kernels \(K_S\) from just the observational measure \(\PP\). Sources describe the most fundamental case in which this is possible, namely, when the causal kernels directly coincide (almost surely) with the corresponding conditional probability measure derived from \(\PP\). 

It was proved \citep[Theorem D.2]{park2023measure} that when one intervenes on \(\sH_U\), then \(\sH_U\) becomes a source. This is a fundamental idea in causality, that when one is able to intervene, then the causal effect of \(\sH_U\) can be obtained by first intervening on \(\sH_U\), and then considering the conditional distribution given \(\sH_U\). We recall this theorem here, since it is used in \ref{sec:proofs} for the proofs of some results in this paper. 
\begin{savetheorem}{thm:fundamental}
    Let \((\Omega,\sH,\PP,\KK)\) be a causal space. Further, let us take an intervention on \(\sH_U\) via \(\QQ\) as in \Cref{def:interventions_causal}, yielding the intervention causal space \((\Omega,\sH,\iPP{U}{\QQ},\iKK{U}{\QQ})\). Then the intervention causal kernel \(\iK{U}{\QQ}{U}\) in the new causal space satisfies the following. 
    \begin{enumerate}[(i)]
        \item It is the same as the corresponding causal \(K_U\) in the original causal space \((\Omega,\sH,\PP,\KK)\) before intervention, i.e. we have
        \begin{equation*}
            K_U=\iK{U}{\QQ}{U};
        \end{equation*}
        \item the causal kernel \(K_U=\iK{U}{\QQ}{U}\) is a version of \(\iPP{U}{\QQ}_{\sH_U}\), which means that \(\sH_U\) is a global source of the intervention causal space. 
    \end{enumerate}
\end{savetheorem}

\section{Proofs}\label{sec:proofs}
\printprop{prop:no_cross_world_conditional}
\begin{proof}
    Let us take any \(S\in\sP(T)\) and any \(\omega\in\Omega\). Suppose that \(\iPP{S}{\delta_\omega}(G)>0\) and \(\iPP{S\setminus U}{\delta_\omega}(G)>0\). Then since \(G\cap A\) and \(A\) both belong to \(\cfact{\sH}\) and \(\sH_U\) has no causal effect on \(\cfact{\sH}\) by the no cross-world causal effect axiom (\Cref{def:counterfactual_causal_space}(ii)), 
    \begin{alignat*}{2}
        \iPP{S}{\delta_\omega}_G(A)&=\frac{K_S(\omega,G\cap A)}{K_S(\omega,G)}\\
        &=\frac{K_{S\setminus U}(\omega,G\cap A)}{K_{S\setminus U}(\omega,G)}\\
        &=\iPP{S\setminus U}{\delta_\omega}_G(A),
    \end{alignat*}
    as required. 
\end{proof}

\printprop{prop:independence_preserved}
\begin{proof}
    \leavevmode\par
    \begin{enumerate}[label=(\roman*)]
        \item See that
        \begin{alignat*}{2}
            \iPP{U}{\QQ}_{\sH_U}(\omega,A\cap B)&\stackrel{a.s.}{=}K_U(\omega,A\cap B)\\
            &\stackrel{(a)}{=}K_U(\omega,A)K_U(\omega,B)\\
            &\stackrel{a.s.}{=}\iPP{U}{\QQ}_{\sH_U}(\omega,A)\iPP{U}{\QQ}_{\sH_U}(\omega,B)
        \end{alignat*}
        where \(\stackrel{a.s.}{=}\) are almost sure equalities that follow from \Cref{thm:fundamental}, and (a) follows from the causal independence of \(A\) and \(B\).
        \item Take arbitrary events \(A\in\fact{\sH}\) and \(B\in\cfact{\sH}\). Then see that
        \begin{alignat*}{2}
            &\iPP{U}{\QQ}(A\cap B)\\
            &\stackrel{(a)}{=}\int\QQ(d\omega_U)K_U(\omega_U,A\cap B)\\
            &\stackrel{(b)}{=}\int\QQ(d\omega_U)K_U(\omega_U,A)K_U(\omega_U,B)\\
            &\stackrel{(c)}{=}\int\QQ((d\omega_{U\cap \fact{T}},d\omega_{U\cap \cfact{T}}))K_{U\cap \fact{T}}(\omega_{U\cap \fact{T}},A)K_{U\cap \cfact{T}}(\omega_{U\cap \cfact{T}},B)\\
            &\stackrel{(d)}{=}\int\QQ(d\omega_{U\cap \fact{T}})K_{U\cap \fact{T}}(\omega_{U\cap \fact{T}},A)\int\QQ(d\omega_{U\cap \cfact{T}})K_{U\cap \cfact{T}}(\omega_{U\cap \cfact{T}},B)\\
            &\stackrel{(e)}{=}\int\QQ(d\omega_U)K_U(\omega_U,A)\int\QQ(d\omega_U)K_U(\omega_U,B)\\
            &=\iPP{U}{\QQ}(A)\iPP{U}{\QQ}(B),
        \end{alignat*}
        where (a) is the definition of the intervention measure, (b) follows from the fact that the factual and counterfactual worlds are causally independent on \(\sH_U\), (c) follows from the no cross-world causal effect axiom (\Cref{def:counterfactual_causal_space}\ref{item:no_cross_world_effect}), (d) follows from the independence of \(\fact{\sH}\) and \(\cfact{\sH}\) under \(\QQ\) and (e) again follows from \Cref{def:counterfactual_causal_space}\ref{item:no_cross_world_effect}. 
    \end{enumerate}
\end{proof}
\printprop{prop:synchronisation_preserved}
\begin{proof}
    See that, by \Cref{thm:fundamental}, for \(\iPP{U}{\QQ}_{\sH_U}\)-almost all \(\omega\in\Omega\),
    \begin{equation*}
        \iPP{U}{\QQ}_{\sH_U}(\omega,A\Delta B)=K_U(\omega,A\Delta B)=0.
    \end{equation*}
    Also, since \(K_U(\omega,A\Delta B)=0\), we have
    \begin{equation*}
        \iPP{U}{\QQ}(A\Delta B)=\int\QQ(d\omega)K_U(\omega,A\Delta B)=0,
    \end{equation*}
    as required. 
\end{proof}

\printtheorem{thm:Nway_interventions}
\begin{proof}
    We check the three axioms of \(N\)-way counterfactual causal spaces given in \Cref{def:Nway_counterfactual_causal_space}\ref{item:Nway_no_cross_world_effect}. 
    \begin{enumerate}[(i)]
        \item For all \(A\in\mathscr{H}\) and \(\omega\in\Omega\), we have, from \Cref{def:intervention_Nway_counterfactual}, 
		\begin{alignat*}{2}
			\iK{U}{\QQ}{\emptyset}(\omega,A)&=\int\QQ(\omega'_U)K_U((\omega_\emptyset,\omega'_U),A)\\
			&=\int\mathbb{Q}(d\omega')K_U(\omega',A)\\
			&=\mathbb{P}^{\text{do}(U,\mathbb{Q})}(A),
		\end{alignat*}
		as required. 
        \item Take an arbitrary \(j\in\{1,...,N\}\), any \(\omega\in\Omega\), any \(S\in\sP(T)\) and any \(A\in\sH^j\). Then we have
        \begin{alignat*}{2}
            \iK{U}{\QQ}{S}(\omega,A)&=\int\QQ(d\omega'_{U\setminus S})K_{S\cup U}((\omega_S,\omega'_{U\setminus S}),A)\\
            &\stackrel{(a)}{=}\int\QQ(d\omega'_{(U\setminus S)\cap T^j})K_{(S\cup U)\cap T^j}((\omega_{S\cap T^j},\omega'_{(U\setminus S)\cap T^j}),A)\\
            &=\int\QQ(d\omega'_{U \setminus(S\cap T^j)}))K_{(S\cap T^j)\cup U}((\omega_{S\cap T^j},\omega'_{U\setminus(S\cap T^j)}),A)\\
            &=\iK{U}{\QQ}{S\cap T^j}(\omega,A)
        \end{alignat*}
        where, in (a), we applied the no cross-world causal effect axiom (\Cref{def:Nway_counterfactual_causal_space}\ref{item:Nway_no_cross_world_effect}). 
        \item For all \(A\in\mathscr{H}_S\) and \(B\in\mathscr{H}\), we have, using the fact that \(A\in\mathscr{H}_S\subseteq\mathscr{H}_{S\cup U}\),
		\begin{alignat*}{2}
			\iK{U}{\QQ}{S}(\omega_S,A\cap B)&=\int\QQ(d\omega'_{U\setminus S})K_{S\cup U}((\omega_S,\omega'_{U\setminus S}),A\cap B)\\
			&=\int\QQ(d\omega'_{U\setminus S})1_A((\omega_S,\omega'_{U\setminus S}))K_{S\cup U}((\omega_S,\omega'_{U\setminus S}),B)\\
			&=\int\QQ(d\omega'_{U\setminus S})1_A(\omega_S)K_{S\cup U}((\omega_S,\omega'_{U\setminus S}),B),
		\end{alignat*}
        where the last line follows because \(1_A\) does not depend on the \(\omega'_{U\setminus S}\) component, as \(A\in\sH_S\). After we take the indicator out of the integration, we are left with
        \[\iK{U}{\QQ}{S}(\omega,A\cap B)=\ind_A(\omega)\int\QQ(d\omega'_{U\setminus S})K_{S\cup U}((\omega_S,\omega'_{U\setminus S}),B).\]
        The integral on the left-hand side is the definition of \(\iK{U}{\QQ}{S}(\omega,B)\). Hence, we have
        \[\iK{U}{\QQ}{S}(\omega,A\cap B)=\ind_A(\omega)\iK{U}{\QQ}{S}(\omega,B),\]
        which is precisely the interventional determinism axiom. 
    \end{enumerate}
\end{proof}

\end{appendix}

\end{document}